\documentclass[journal]{IEEEtran}
\usepackage{cite}

\ifCLASSINFOpdf
  \usepackage[pdftex]{graphicx}
  \usepackage{epstopdf}
  \usepackage{subfigure}
  \graphicspath{{./}{./figures/}}
\else
\fi
\usepackage{xcolor}

\usepackage{amsmath}
\usepackage{amssymb}

\newcommand*{\qed}{\hfill\ensuremath{\square}}%
\newcommand*\diff{\mathop{}\!\mathrm{d}}%

\newtheorem{theorem}{Theorem}[section]
\newtheorem{corollary}{Corollary}[theorem]
\newtheorem{lemma}[theorem]{Lemma}

\newtheorem{remark}[theorem]{Remark}
\newtheorem{definition}[theorem]{Definition}

\allowdisplaybreaks

\begin{document}
%
\title{Boundary input-to-state stabilization of a damped Euler-Bernoulli beam in the presence of a state-delay}
%
%
%

\author{Hugo~Lhachemi, Robert Shorten
\thanks{This publication has emanated from research supported in part by a research grant from Science Foundation Ireland (SFI) under grant number 16/RC/3872 and is co-funded under the European Regional Development Fund and by I-Form industry partners.
\newline\indent Hugo Lhachemi is with the School of Electrical and Electronic Engineering, University College Dublin, Dublin, Ireland (e-mail: hugo.lhachemi@ucd.ie).
\newline\indent Robert Shorten is with the School of Electrical and Electronic Engineering, University College Dublin, Dublin, Ireland, and also with the Dyson School of Design Engineering, Imperial College London, London, U.K. (e-mail: robert.shorten@ucd.ie).
}
}

\maketitle

\begin{abstract}
This paper is concerned with the point torque boundary feedback stabilization of a damped Euler-Bernoulli beam model in the presence of a time-varying state-delay. First, a finite-dimensional truncated model is derived by spectral reduction. Then, for a given stabilizing state-feedback of the delay-free truncated model, an LMI-based sufficient condition on the maximum amplitude of the delay is employed to guarantee the stability of the closed-loop state-delayed truncated model. Second, we assess the exponential stability of the resulting closed-loop infinite-dimensional system under the assumption that the number of modes of the original infinite-dimensional system captured by the truncated model has been selected large enough. Finally, we consider in our control design the possible presence of a distributed perturbation, as well as additive boundary perturbations in the control inputs. In this case, we derive for the closed-loop system an exponential input-to-state estimate with fading memory of the distributed and boundary disturbances.
\end{abstract}


%
\IEEEpeerreviewmaketitle

\section{Introduction}\label{sec: Introduction}
%
%
%
%

The study of the stability properties of various partial differential equations (PDEs) models under either boundary or state delays is an active topic of research~\cite{nicaise2008stabilization,nicaise2007stabilization,nicaise2009stability,fridman2009exponential,solomon2015stability}.  In this context, the stabilization of open-loop unstable PDEs in the presence of a delayed term has attracted much attention in recent years. In this context, one of the early works reported in the literature deals with the boundary feedback stabilization of an unstable reaction-diffusion equation under large constant input delays by means of a backstepping transformation~\cite{krstic2009control}. More recently, a similar problem was tackled in~\cite{prieur2018feedback} by resorting to finite-dimensional controller synthesis methods. First, a finite-dimensional model capturing the unstable modes of the original plant was derived via spectral decomposition. Then, the controller was obtained by resorting to a predictor feedback~\cite{krstic2009delay}, which is an efficient tool for the feedback stabilization of finite-dimensional linear time-invariant systems with constant input delays. Finally, the stability of the resulting closed-loop infinite-dimensional system was assessed via the use of an adequate Lyapunov functional. This control strategy was reused in~\cite{guzman2019stabilization} for the delay boundary stabilization of a linear Kuramoto-Sivashinsky equation, and then generalized in~\cite{lhachemi2019feedback,lhachemi2019control,lhachemi2019lmi} to a class of diagonal infinite-dimensional systems for either constant or time-varying input delays.

Going beyond the delay boundary stabilization of PDEs, the boundary stabilization of PDEs in the presence of a state-delay has also been the subject of a number of recent publications. In this context, most of the reported works deal with the boundary feedback stabilization of an unstable reaction-diffusion equation in the presence of a state-delay in the reaction term. The case of a constant state-delay with Dirichlet boundary conditions was reported in~\cite{hashimoto2016stabilization} where the control design was performed via a backstepping transformation. The case of a time-varying state-delay was investigated in~\cite{kang2017boundaryIFAC} for Neumann-Dirichlet boundary conditions with Dirichlet actuation. The control design also took advantage of a backstepping transformation while resorting to an LMI-based argument for assessing the stability of the closed-loop system. This work was later extended in~\cite{kang2017boundary} for the stabilization of a cascade PDE-ODE system with either Dirichlet or Neumann actuation. The case of Robin boundary conditions was studied in~\cite{lhachemi2019boundary} in the presence of both a time-varying state-delay and a distributed disturbance. Following in a similar manner to the idea reported in~\cite{prieur2018feedback} for the delay boundary stabilization of a reaction-diffusion equation, the proposed control strategy consists in designing the controller on a finite-dimensional truncated model capturing a finite number of modes of the original plant. The subsequent stability analysis showed that the resulting closed-loop system is input-to-state stable (ISS)~\cite{sontag1989smooth} with fading memory~\cite{karafyllis2019input} of the distributed perturbation. It is worth noting that very few works have been reported on the extension of the aforementioned methods to other types of boundary control systems presenting a state-delay. Among the reported works in this field, one can find the case of a linear~\cite{kang2018regional} and nonlinear~\cite{kang2018boundary} Schr{\"o}dinger equation by means of a backstepping transformation and an LMI-based stability analysis. Extensions to other types of PDEs, such as wave and beam equations, remain open.

In this paper, we are concerned with the point torque boundary feedback stabilization of a damped Euler-Bernoulli beam in the presence of both a time-varying state-delay and additive perturbations. Specifically, we aim at achieving the boundary input-to-state stabilization of the system with respect to both distributed and boundary disturbances. The main motivation of such a goal relies in the fact that the ISS property, originally introduced by Sontag in~\cite{sontag1989smooth}, is one of the main tools for assessing the robustness of a system with respect to disturbances. This property can also be used for the establishment of small gain conditions ensuring the stability of interconnected systems~\cite{karafyllis2019input}. The establishment of ISS properties for finite-dimensional systems was intensively studied during the last three decades. However, its extension to infinite-dimensional systems is more recent, particularly in the case of boundary perturbations, and remains highly challenging~\cite{argomedo2012d,jacob2018infinite,jacob2018continuity,
karafyllis2016iss,karafyllis2017iss,karafyllis2019input,
lhachemi2019input,lhachemi2018iss,
mironchenko2019monotonicity,mironchenko2016restatements,
mironchenko2017characterizations,zheng2017giorgi,zheng2017input}. 
It is also worth noting that most of these works deal with the establishment of ISS properties with respect to boundary disturbances for open loop stable systems. The literature regarding the feedback input-to-state stabilization of open loop unstable infinite-dimensional systems with respect to boundary perturbations is less developed~\cite{tanwani2017disturbance}.

In this context, the present paper is concerned with the design of a point torque boundary feedback control law that ensures the exponential input-to-state stabilization, with respect to both distributed and boundary perturbations, of a damped Euler-Bernoulli beam presenting a time-varying state-delay. The adopted approach, which allows either one single command input (located at one of the two boundaries of the domain) or two command inputs, is organized as follows. First, a spectral decomposition of the studied beam model is carried out. This spectral decomposition is used for deriving a finite-dimensional truncated model of the beam capturing its unstable modes plus an adequately chosen number of slow stable modes. Specifically, by means of a small gain argument, the order of the truncated model is selected in order to ensure the robust stability of the residual infinite-dimensional system with respect to exponentially vanishing command inputs exhibiting a prescribed decay rate. In this context, the proposed control law consists in a state feedback of the truncated model. The stability of the resulting closed-loop truncated model in the presence of a state-delay is assessed via an LMI-based (sufficient) condition on the maximum amplitude of the state-delay. Under the assumption of an adequate choice of the order of the truncated model, we show that the resulting infinite-dimensional closed-loop system satisfies an exponential ISS-like estimate\footnote{We refer the reader to Remark~\ref{rem: ISS-like estimate} for the meaning of ``ISS-like''.} with fading memory of the perturbations.

This paper is organized as follows. The problem setting and the resulting abstract boundary control system are introduced in Section~\ref{eq: problem setting and abstract form}. After introducing the proposed control strategy, the main result of this paper is stated in Section~\ref{sec: control strat and main result}. The subsequent stability analysis is carried out in two steps: first in Section~\ref{sec: stab CL finite-dim} for the stability analysis of the truncated model and then in Section~\ref{sec: stab CL infinite-dim} for the stability analysis of the resulting closed-loop infinite-dimensional system. The obtained results are numerically illustrated in Section~\ref{sec: numerical illustration}. Finally, concluding remarks are provided in Section~\ref{sec: conclusion}.

\section{Problem setting and abstract form}\label{eq: problem setting and abstract form}

\subsection{Notations and definitions}

The sets of non-negative integers, positive integers, real, non-negative real, positive real, and complex numbers are denoted by $\mathbb{N}$, $\mathbb{N}^*$, $\mathbb{R}$, $\mathbb{R}_+$, $\mathbb{R}_+^*$, and $\mathbb{C}$, respectively. The real and imaginary parts of a complex number $z$ are denoted by $\operatorname{Re} z$ and $\operatorname{Im} z$, respectively. The field $\mathbb{K}$ denotes either $\mathbb{R}$ or $\mathbb{C}$. The set of $n$-dimensional vectors over $\mathbb{K}$ is denoted by $\mathbb{K}^n$ and is endowed with the Euclidean norm $\Vert x \Vert = \sqrt{x^* x}$. The set of $n \times m$ matrices over $\mathbb{K}$ is denoted by $\mathbb{K}^{n \times m}$ and is endowed with the induced norm denoted by $\Vert\cdot\Vert$. For any symmetric matrix $P \in \mathbb{R}^{n \times n}$, $P \succ 0$ (resp. $P \succeq 0$) means that $P$ is positive definite (resp. positive semi-definite). The set of symmetric positive definite matrices of order $n$ is denoted by $\mathbb{S}_n^{+*}$. For any symmetric matrix $P \in \mathbb{R}^{n \times n}$, $\lambda_m(P)$ and $\lambda_M(P)$ denote the smallest and largest eigenvalues of $P$, respectively.

The set of square-integrable functions (w.r.t. the Lebesgue measure) over the interval $(0,1) \subset \mathbb{R}$ is denoted by $L^2(0,1)$ and is endowed with its natural inner product $\left[ f , g \right] = \int_{0}^{1} f(\xi) \overline{g(\xi)} \diff\xi$, providing a structure of Hilbert space. Denoting by $f'$, when it exists, the weak derivative of $f \in L^2(0,1)$, we consider the Sobolev space $H^m(0,1) \triangleq \{ f \in L^2(0,1) \, : \, f',f'', \ldots , f^{(m)} \in L^2(0,1) \}$. Finally, we introduce $H_0^1(0,1) \triangleq \{ f \in H^1(0,1) \; : \; f(0)=f(1)=0 \}$.

\subsection{Problem setting and control objective}

As we have mentioned, our focus is the point torque boundary feedback stabilization of the following damped Euler Bernoulli beam model in the presence of a state-delay:
\begin{subequations}
\begin{align}
y_{tt}(t,x) & = - y_{xxxx}(t,x) + 2 \alpha y_{txx}(t,x) \label{eq: beam model - PDE} \\
& \phantom{=}\; + \beta_0 y(t,x) + \gamma y(t-h(t),x) + d_d(t,x) \\
y(t,0) & = y(t,1) = 0 \label{eq: beam model - zero BC} \\
y_{xx} (t,0) & = u_1(t) + d_{b,1}(t) \label{eq: beam model - left torque input} \\
y_{xx} (t,1) & = u_2(t) + d_{b,2}(t) \label{eq: beam model - right torque input} \\
y(\tau,x) & = y_0(\tau,x) , \;\; \tau\in[-h_M,0] \label{eq: beam model - IC1} \\
y_{t} (\tau,x) & = y_{t0}(\tau,x) , \;\; \tau\in[-h_M,0] \label{eq: beam model - IC2}
\end{align}
\end{subequations}
for $t > 0$ and $x \in (0,1)$. Here we have the damping parameter $\alpha > 1$ and the coefficients $\beta_0 \geq 0$ and $\gamma > 0$. We also introduce the quantity $\beta = \beta_0 + \gamma > 0$ that will be used in the sequel. The boundary conditions (\ref{eq: beam model - left torque input}-\ref{eq: beam model - right torque input}) involve the point torque boundary control inputs $u_1,u_2 : \mathbb{R}_+ \rightarrow \mathbb{R}$ (with possibly one identically equal to zero), $h: \mathbb{R}_+ \rightarrow \mathbb{R}$ is a time-varying delay, $d_d : \mathbb{R}_+ \times (0,1) \rightarrow \mathbb{R}$ is a distributed disturbance, and $d_{b,1},d_{b,2} : \mathbb{R}_+ \rightarrow \mathbb{R}$ are boundary disturbances. We assume that there exist constants, $0 < h_m < h_M$, such that $h_m \leq h(t) \leq h_M$ for all $t \geq 0$. Finally, $y_0,y_{t0}: [-h_M,0] \times (0,1) \rightarrow \mathbb{R}$ are the initial conditions.

The tackled control objective consists of designing the boundary control inputs $u_1,u_2$ in order to ensure the closed-loop stability of the Euler-Bernoulli beam (\ref{eq: beam model - PDE}-\ref{eq: beam model - IC2}) for any continuous time varying delay such that $0 < h_m \leq h \leq h_M$ with $h_M > 0$ to be characterized. Moreover, this control strategy must ensure the ISS property of the closed-loop system with respect to both distributed and boundary perturbations.

\subsection{Abstract form}

We consider the state-space of the form of the Hilbert space $\mathcal{H} = \left( H^2(0,1) \cap H_0^1(0,1) \right) \times L^2(0,1)$ with associated inner product defined for all $(y_1,y_2),(\hat{y}_1,\hat{y}_2)\in\mathcal{H}$ by 
\begin{align*}
\left<(y_1,y_2) , (\hat{y}_1,\hat{y}_2) \right>
& = \int_0^1 y_1''(x) \overline{\hat{y}_1''(x)} + y_2(x) \overline{\hat{y}_2(x)} \diff x .
\end{align*}
We also introduce the bounded map $\Pi \in \mathcal{L}(\mathcal{H})$ defined by $\Pi (y_1,y_2) = (0,y_1)$. 

\begin{remark}\label{rem: Pi bounded}
The bounded nature of $\Pi$ follows from the fact that, as $y_1 \in H^2(0,1) \cap H_0^1(0,1)$, it is continuously differentiable with $y_1(0)=y_1(1)=0$. Rolle's theorem provides $a \in (0,1)$ such that $y'_1(a)=0$. Then, for all $x \in [0,1]$, we have $y_1(x) = \int_0^x \int_a^\xi y''_1(s) \diff s \diff \xi$. We infer by Cauchy-Schwarz $\vert y_1(x) \vert^2 \leq \int_0^1 \vert y''_1(s) \vert^2 \diff s$. This yields $\Vert \Pi (y_1,y_2) \Vert = \sqrt{\int_0^1 \vert y_1(x) \vert^2 \diff x} \leq \Vert (y_1,y_2) \Vert$.
\end{remark}

The beam model (\ref{eq: beam model - PDE}-\ref{eq: beam model - IC2}) can be written as the abstract boundary control system: 
\begin{subequations}
\begin{align}
\dfrac{\mathrm{d}X}{\mathrm{d}t}(t) & = \mathcal{A}X(t) + \gamma \Pi X(t-h(t)) + p_d(t) \label{eq: abstract form - DE} \\
\mathcal{B}X(t) & = w(t) \triangleq u(t) + d_b(t) \label{eq: abstract form - BC} \\
X(\tau) & = \Phi(\tau), \quad \tau \in [-h_M,0] \label{eq: abstract form - IC} 
\end{align}
\end{subequations}
for $t>0$ where $\mathcal{A}(y_1,y_2) = ( y_2 , - y''''_1 + 2 \alpha y''_2 + \beta_0 y_1 ) \in\mathcal{H}$ defined on $D(\mathcal{A}) = 
\left( H^4(0,1) \cap H_0^1(0,1) \right) \times \left( H^2(0,1) \cap H_0^1(0,1) \right)$ and $\mathcal{B}(y_1,y_2) = (y''_1(0),y''_1(1)) \in \mathbb{R}^2$ defined on $D(\mathcal{B}) = D(\mathcal{A})$. Here we have the state $X(t) = (y(t,\cdot),y_t(t,\cdot)) \in \mathcal{H}$, the control input $u(t)=(u_1(t),u_2(t)) \in \mathbb{R}^2$, the distributed disturbance $p_d(t) = (0,d_d(t,\cdot)) \in \mathcal{H}$, the bounday perturbation $d_b(t) = (d_{b,1}(t),d_{b,2}(t)) \in \mathbb{R}^2$, and the initial condition $\Phi(t) = (y_0(t,\cdot),y_{t0}(t,\cdot)) \in \mathcal{H}$.

Following the terminology of~\cite[Def.~3.3.2]{Curtain2012}, $(\mathcal{A},\mathcal{B})$ is an abstract boundary control system. Indeed, the disturbance-free operator $\mathcal{A}_0 \triangleq \left.\mathcal{A}\right\vert_{D(\mathcal{A}_0)}$ with $D(\mathcal{A}_0) \triangleq D(\mathcal{A}) \cap \mathrm{ker}(\mathcal{B})$ generates a $C_0$-semigroup $S(t)$. Furthermore, introducing $L:\mathbb{R}^2 \rightarrow \mathcal{H}$ defined for any $u = (u_1,u_2) \in \mathbb{R}^2$ by $[Lu](x) = \left( \dfrac{u_2-u_1}{6}x^3+\dfrac{u_1}{2}x^2-\dfrac{2 u_1 + u_2}{6}x , 0 \right)$ for any $x \in [0,1]$, $L \in \mathcal{L}(\mathbb{R}^2,\mathcal{H})$ is a lifting operator in the sense that $R(L) \subset D(\mathcal{A})$, $\mathcal{A} L$ is bounded, and $\mathcal{B}L = I_{\mathbb{R}^2}$.

In this paper, we consider the concept of mild solutions of the abstract boundary control problem (\ref{eq: abstract form - DE}-\ref{eq: abstract form - IC}). Assuming that $p_d \in \mathcal{C}^0(\mathbb{R}_+;\mathcal{H})$, $u , d_b \in \mathcal{C}^1(\mathbb{R}_+;\mathbb{R}^2)$, $h \in \mathcal{C}^0(\mathbb{R}_+;\mathbb{R})$ with $0 < h_m \leq h \leq h_M$, and $\Phi \in \mathcal{C}^0([-h_M,0];\mathcal{H})$, the mild solution $X \in \mathcal{C}^0(\mathbb{R}_+;\mathcal{H})$ of (\ref{eq: abstract form - DE}-\ref{eq: abstract form - IC}) is uniquely defined by
\begin{align}
& X(t) 
= S(t) \{ \Phi(0) - L w(0) \} + L w(t) \label{eq: def mild solution} \\
& \, + \int_0^t S(t-s) \{ \mathcal{A}Lw(s) - L\dot{w}(s) + \gamma \Pi X(s-h(s)) + p_d(s) \} \diff s \nonumber
\end{align}
for $t \geq 0$ with $w = u+d_b$ and the initial condition $X(\tau) = \Phi(\tau)$ for all $\tau\in[-h_M,0]$.

\subsection{Spectral properties of the beam}

In preparation for control design, we rewrite (\ref{eq: abstract form - DE}) under the following form:
\begin{equation}\label{eq: abstract form - DE bis}
\dfrac{\mathrm{d}X}{\mathrm{d}t}(t) = \mathcal{U} X(t) + \gamma \Pi \left\{ X(t-h(t)) - X(t) \right\} + p_d(t) ,
\end{equation} 
where $\mathcal{U} \triangleq \mathcal{A} + \gamma \Pi$ with $D(\mathcal{U}) = D(\mathcal{A})$, i.e., introducing $\beta = \beta_0 + \gamma > 0$, $\mathcal{U}(y_1,y_2) = ( y_2 , - y''''_1 + 2 \alpha y''_2 + \beta y_1 )$. We also introduce $\mathcal{U}_0 = \mathcal{A}_0 + \gamma \Pi$ defined on $D(\mathcal{U}_0) = D(\mathcal{A}_0)$ the disturbance-free operator associated with $(\mathcal{U},\mathcal{B})$ and $T(t)$ the $C_0$-semigroup generated by $\mathcal{U}_0$. Then the mild solution (\ref{eq: def mild solution}) can be equivalently rewritten in function of $T$ as (see Appendix~\ref{annex: mild solution}):
\begin{align}
X(t) 
& = T(t) \{ \Phi(0) - L w(0) \} + L w(t) \label{eq: def mild solution bis} \\
& \phantom{=}\; + \int_0^t T(t-s) \bigg\{ \mathcal{U}Lw(s) - L\dot{w}(s) \nonumber \\
& \hspace{2cm} + \gamma \Pi \{ X(s-h(s)) - X(s) \} + p_d(s) \bigg\} \diff s \nonumber
\end{align}
for $t \geq 0$.

The eigenstructures of $\mathcal{U}_0$ are characterized by the following lemma. The proof is a straightforward extension to the case $\beta > 0$ of the approach reported in~\cite[Exercise~2.23]{Curtain2012} for $\beta = 0$.

\begin{lemma}
The point spectrum of $\mathcal{U}_0$ is given by $\sigma_p(\mathcal{U}_0) = \left\{ \lambda_{n,\epsilon} \,:\, n \in \mathbb{N}^* , \, \epsilon \in \{-1,+1\} \right\}$ with simple eigenvalues
\begin{equation}\label{eq: eigenvalue U0}
\lambda_{n,\epsilon} = - \alpha n^2 \pi^2 + \epsilon \sqrt{(\alpha^2-1)n^4\pi^4+\beta} \in \mathbb{R} 
\end{equation}
and associated unit eigenvectors
\begin{equation*}
\phi_{n,\epsilon} = \dfrac{1}{k_{n,\epsilon}} \left( \sin(n\pi\cdot) , \lambda_{n,\epsilon} \sin(n\pi\cdot) \right)
\end{equation*}
where $k_{n,\epsilon} > 0$ is given by $k_{n,\epsilon} = \sqrt{(n^4 \pi^4 + \lambda_{n,\epsilon}^2)/2}$.
\end{lemma}

\begin{remark}\label{rem: properties eigenvalues}
From (\ref{eq: eigenvalue U0}), it is easy to see that the following hold:
\begin{itemize}
\item $\lambda_{n,\epsilon}^2 + 2\alpha n^2 \pi^2 \lambda_{n,\epsilon} + (n^4 \pi^4 - \beta) = 0$ for all $n \geq 1$ and $\epsilon \in \{ -1 , +1 \}$;
\item $\lambda_{n,\epsilon} < \lambda_{m,\epsilon}$ for all $n > m \geq 1$ and $\epsilon \in \{ -1 , +1 \}$;
\item $\lambda_{n,-1} < 0$ and $\lambda_{n,-1} < \lambda_{n,+1}$ for all $n \geq 1$;
\item $\lambda_{n,+1} \geq 0$ if and only if $1 \leq n \leq \beta^{1/4}/\pi$.
\end{itemize} 
In particular, the number of unstable eigenvalues of the disturbance-free operator $\mathcal{U}_0$ is given by $\lfloor \beta^{1/4}/\pi \rfloor$.
\end{remark}

One of the key properties that will be used in the sequel is the concept of Riesz basis~\cite{christensen2016introduction} which is recalled in the following definition.

\begin{definition}[Riesz basis~\cite{christensen2016introduction}]\label{def: Riesz basis}
A sequence $\mathcal{F} = \{ \varphi_k ,\; k\in\mathbb{N}\}$ of vectors of a Hilbert space $X$ over $\mathbb{K}$ is a Riesz basis if 
\begin{enumerate}
\item $\mathcal{F}$ is maximal: $\overline{\mathrm{span}_\mathbb{K}(\mathcal{F})} = X$, i.e., the closure of the vector space spanned by $\mathcal{F}$ coincides with the whole space $X$;
\item there exist $m_R,M_R\in\mathbb{R}_+^*$ such that for any $N\in\mathbb{N}$ and any $a_k \in\mathbb{K}$, 
\end{enumerate}
\begin{equation}\label{eq: def Riesz basis inequality}
m_R \sum\limits_{0 \leq k \leq N} |a_k|^2
\leq
\left\Vert \sum\limits_{0 \leq k \leq N} a_{k} \varphi_{k} \right\Vert^2
\leq
M_R \sum\limits_{0 \leq k \leq N} |a_k|^2 . 
\end{equation}
\end{definition}

We can now introduce the following lemma, which will be crucial in the sequel and whose proof is placed in Appendix~\ref{annex: Riesz basis}.
\begin{lemma}\label{lem: Riesz basis}
$\mathcal{F}_\phi = \left\{ \phi_{n,\epsilon} \,:\, n \in \mathbb{N}^* , \, \epsilon \in \{-1,+1\} \right\}$ is a Riesz basis of $\mathcal{H}$.
\end{lemma}

In particular, following the therminology of~\cite[Def.~2.3.4]{Curtain2012}, $\mathcal{U}_0$ is a Riesz-spectral operator. We introduce $\mathcal{F}_\psi = \left\{ \psi_{n,\epsilon} \,:\, n \in \mathbb{N}^* , \, \epsilon \in \{-1,+1\} \right\}$ the dual Riesz basis of $\mathcal{F}_\phi$: 
\begin{equation*}
\psi_{n,\epsilon} 
= C_{n,\epsilon} \left( -\dfrac{\lambda_{n,-\epsilon}}{n^4\pi^4} \sin(n\pi\cdot) , \sin(n\pi\cdot) \right) 
\end{equation*}
with
\begin{equation*}
C_{n,\epsilon}
= \dfrac{2 k_{n,\epsilon}}{\lambda_{n,\epsilon}-\lambda_{n,-\epsilon}}
= \dfrac{\epsilon k_{n,\epsilon}}{\sqrt{(\alpha^2-1)n^4\pi^4+\beta}} 
\neq 0 ,
\end{equation*}
i.e., $\left< \phi_{n_1,\epsilon_1} , \psi_{n_2,\epsilon_2} \right> = \delta_{(n_1,\epsilon_1),(n_2,\epsilon_2)} \in \{ 0,1 \}$ with $\delta_{(n_1,\epsilon_1),(n_2,\epsilon_2)} = 1$ if and only if $(n_1,\epsilon_1)=(n_2,\epsilon_2)$. From the general theory of Riesz bases (see, e.g., \cite{christensen2016introduction}), we have
\begin{equation*}
\forall z \in \mathcal{H} , \;
z = \sum\limits_{\substack{n \geq 1 \\ \epsilon\in\{-1,+1\}}} \left< z , \psi_{n,\epsilon} \right> \phi_{n,\epsilon} 
\end{equation*}
with
\begin{equation}\label{eq: Riesz basis inequality}
m_R \sum\limits_{\substack{n \geq 1 \\ \epsilon\in\{-1,+1\}}} \vert \left< z , \psi_{n,\epsilon} \right> \vert^2
\leq \Vert z \Vert^2
\leq M_R \sum\limits_{\substack{n \geq 1 \\ \epsilon\in\{-1,+1\}}} \vert \left< z , \psi_{n,\epsilon} \right> \vert^2 
\end{equation}
where the constants $m_r,M_R > 0$ are provided in Appendix~\ref{annex: Riesz basis}. Furthermore, as $\mathcal{U}_0$ is a Riesz-spectral operator generating the $C_0$-semigroup $T(t)$, we have~\cite[Thm.~2.3.5]{Curtain2012}
\begin{equation}\label{eq: expression C0-semigroup T}
\forall z \in \mathcal{H} , \; \forall t \geq 0 , \;
T(t)z = \sum\limits_{\substack{n \geq 1 \\ \epsilon\in\{-1,+1\}}} e^{\lambda_{n,\epsilon}t} \left< z , \psi_{n,\epsilon} \right> \phi_{n,\epsilon}
\end{equation}

\section{Spectral decomposition, control strategy, and main result}\label{sec: control strat and main result}

\subsection{Spectral decomposition}\label{subsec: spectral decomposition}
We define $c_{n,\epsilon}(t) = \left< X(t) , \psi_{n,\epsilon} \right>$ the coefficients of projection of the system trajectory $X$ into the Riesz basis $\mathcal{F}_\phi$.

\subsubsection{Preliminary spectral decomposition}
Let $X = (x_1,x_2) \in \mathcal{C}^0(\mathbb{R}_+;\mathcal{H})$ be a mild solution (\ref{eq: def mild solution}) of (\ref{eq: abstract form - DE}-\ref{eq: abstract form - IC}). Then $c_{n,\epsilon}\in\mathcal{C}^0(\mathbb{R}_+;\mathbb{R})$ and, using (\ref{eq: def mild solution bis}), (\ref{eq: expression C0-semigroup T}), and the integration by parts:
\begin{align*}
& \int_0^t e^{\lambda_{n,\epsilon}(t-s)} \left< L \dot{w}(s),\psi_{n,\epsilon} \right> \diff s \\
& \qquad = \left< L w(t),\psi_{n,\epsilon} \right> - e^{\lambda_{n,\epsilon}t} \left< L w(0),\psi_{n,\epsilon} \right> \\
& \qquad \phantom{=}\; + \lambda_{n,\epsilon} \int_0^t e^{\lambda_{n,\epsilon}(t-s)} \left< L w(s),\psi_{n,\epsilon} \right> \diff s ,
\end{align*}
we have for all $t \geq 0$ 
\begin{align*}
c_{n,\epsilon}(t)
& = e^{\lambda_{n,\epsilon}t} c_{n,\epsilon}(0) + \int_0^t e^{\lambda_{n,\epsilon}(t-s)} f(s) \diff s
\end{align*}
with $f(t) = \big< - \lambda_{n,\epsilon}Lw(t) + \mathcal{U}Lw(t) + \gamma\Pi \{ X(t-h(t)) - X(t) \} + p_d(t) , \psi_{n,\epsilon} \big>$. As $f$ is continuous, we have $c_{n,\epsilon}\in\mathcal{C}^1(\mathbb{R}_+;\mathbb{R})$ and the following ODE (see also~\cite{lhachemi2018iss}) is satisfied for all $t \geq 0$
\begin{align}
\dot{c}_{n,\epsilon}(t)
& = \lambda_{n,\epsilon} c_{n,\epsilon}(t) + \gamma \Delta_{n,\epsilon}(t) \label{eq: spectral decomposition - prel}\\
& \phantom{=}\; - \lambda_{n,\epsilon} \left< L \{u(t)+d_b(t)\} , \psi_{n,\epsilon} \right> \nonumber \\
& \phantom{=}\; + \left< \mathcal{U} L \{u(t)+d_b(t)\} , \psi_{n,\epsilon} \right> + p_{d,n,\epsilon}(t) , \nonumber
\end{align}
where, using $\psi_{n,\epsilon} = (\psi_{n,\epsilon}^1 , \psi_{n,\epsilon}^2)$, we have $\Delta_{n,\epsilon}(t) = \big< \Pi \{ X(t-h(t)) - X(t) \} , \psi_{n,\epsilon} \big> = \big[ x_1(t-h(t)) - x_1(t) , \psi_{n,\epsilon}^2 \big]$, and $p_{d,n,\epsilon}(t) = \big< p_d(t) , \psi_{n,\epsilon} \big> = \big[ d_d(t) , \psi_{n,\epsilon}^2 \big]$. 

The main idea of the control design consists in using the spectral decomposition (\ref{eq: spectral decomposition - prel}) to obtain a truncated finite-dimensional model capturing a sufficient number of modes of the original infinite-dimensional system (\ref{eq: beam model - PDE}-\ref{eq: beam model - IC2}). However, the ODE (\ref{eq: spectral decomposition - prel}) is not yet under a suitable form for control design because the term $\Delta_{n,\epsilon}$ depends on the state trajectory $X$. Thus we need to express this term in function of the coefficients of projection $c_{n,\epsilon}$.

\subsubsection{Expression of the term $\Delta_{n,\epsilon}$ in function of the coefficients of projection $c_{n,\epsilon}$}
First, we note that, for all $n \geq 1$,
\begin{equation*}
\psi_{n,-1}^1 = - \dfrac{k_{n,-1}\lambda_{n,+1}}{k_{n,+1}\lambda_{n,-1}} \psi_{n,+1}^1 , 
\quad
\psi_{n,-1}^2 = - \dfrac{k_{n,-1}}{k_{n,+1}} \psi_{n,+1}^2 ,
\end{equation*}
where the above identity are well defined because $k_{n,+1} \neq 0$ and $\lambda_{n,-1} \neq 0$. We obtain that
\begin{equation*}
\begin{bmatrix}
c_{n,+1} \\ c_{n,-1}
\end{bmatrix}
=
\begin{bmatrix}
1 & 1 \\ - \dfrac{k_{n,-1}\lambda_{n,+1}}{k_{n,+1}\lambda_{n,-1}} & - \dfrac{k_{n,-1}}{k_{n,+1}}
\end{bmatrix}
\begin{bmatrix}
\left[ x_1'' , (\psi_{n,+1}^{1})'' \right] \\ \left[ x_2 , \psi_{n,+1}^2 \right]
\end{bmatrix} .
\end{equation*}
The determinant of the above square matrix is given by
\begin{equation*}
\mathrm{det}_n = \dfrac{k_{n,-1}}{k_{n,+1}} \left\{ \dfrac{\lambda_{n,+1}}{\lambda_{n,-1}} - 1 \right\} \neq 0 ,
\end{equation*}
thus
\begin{equation*}
\begin{bmatrix}
\left[ x_1'' , (\psi_{n,+1}^{1})'' \right] \\ \left[ x_2 , \psi_{n,+1}^2 \right]
\end{bmatrix}
=
\dfrac{1}{\mathrm{det}_n}
\begin{bmatrix}
- \dfrac{k_{n,-1}}{k_{n,+1}} & -1 \\ \dfrac{k_{n,-1}\lambda_{n,+1}}{k_{n,+1}\lambda_{n,-1}} & 1
\end{bmatrix}
\begin{bmatrix}
c_{n,+1} \\ c_{n,-1}
\end{bmatrix}
\end{equation*}
and we infer that
\begin{equation*}
\left[ x_1'' , (\psi_{n,+1}^{1})'' \right]
= \dfrac{\lambda_{n,-1}}{k_{n,-1}} \times \dfrac{k_{n,-1} c_{n,+1} + k_{n,+1} c_{n,-1}}{\lambda_{n,-1}-\lambda_{n,+1}} .
\end{equation*}
Now, we note that
\begin{equation*}
\left[ x_1 , \psi_{n,\epsilon}^2 \right]
= C_{n,\epsilon} \left[ x_1 , \sin(n\pi\cdot) \right]
= \dfrac{\epsilon k_{n,\epsilon}}{k_{n,+1}} \left[ x_1 , \psi_{n,+1}^2 \right]
\end{equation*}
and, using two integration by parts, the boundary conditions $x_1(t,0)=x_1(t,1)=0$, and the identity $(\psi_{n,+1}^1)'' = \dfrac{C_{n,+1}\lambda_{n,-1}}{n^2 \pi^2} \sin(n\pi\cdot)$, 
\begin{align*}
\left[ x_1 , \psi_{n,+1}^2 \right]
& = C_{n,+1} \int_{0}^{1} x_1(t,\xi) \sin(n\pi\xi) \diff\xi \\
& = - \dfrac{C_{n,+1}}{n^2 \pi^2} \int_{0}^{1} x''_1(t,\xi) \sin(n\pi\xi) \diff\xi \\
& = - \dfrac{1}{\lambda_{n,-1}} \left[ x''_1 , (\psi_{n,+1}^1)'' \right] .
\end{align*}
Then we have
\begin{align*}
\left[ x_1 , \psi_{n,\epsilon}^2 \right]
& = - \dfrac{\epsilon k_{n,\epsilon}}{\lambda_{n,-1} k_{n,+1}} \left[ x''_1 , (\psi_{n,+1}^1)'' \right] \\
& = - \dfrac{\epsilon k_{n,\epsilon}}{k_{n,-1} k_{n,+1}} \times \dfrac{k_{n,-1} c_{n,+1} + k_{n,+1} c_{n,-1}}{\lambda_{n,-1}-\lambda_{n,+1}} .
\end{align*}
Hence the term $\Delta_{n,\epsilon}(t)$ can be rewritten in function of the coefficients of projection $c_{n,-1}(t)$ and $c_{n,+1}(t)$ as follows:
\begin{align}
\Delta_{n,\epsilon}(t)
& = - \epsilon \dfrac{k_{n,\epsilon}}{k_{n,-1}} \times \dfrac{c_{n,-1}(t-h(t)) - c_{n,-1}(t)}{\lambda_{n,-1}-\lambda_{n,+1}} \label{eq: expression Delta,n,eps} \\
& \phantom{=}\; - \epsilon \dfrac{k_{n,\epsilon}}{k_{n,+1}} \times \dfrac{c_{n,+1}(t-h(t)) - c_{n,+1}(t)}{\lambda_{n,-1}-\lambda_{n,+1}} . \nonumber
\end{align}

\subsubsection{Spectral decomposition in suitable form for control design}
We can now derive from the preliminary spectral decomposition (\ref{eq: spectral decomposition - prel}) and equation (\ref{eq: expression Delta,n,eps}) a spectral decomposition that is suitable for control design. Let $\{e_1,e_2\}$ be the canonical basis of $\mathbb{R}^2$ and define 
\begin{equation}\label{eq: coeff b_n}
b_{n,\epsilon,m} = - \lambda_{n,\epsilon} \left< L e_m , \psi_{n,\epsilon} \right> + \left< \mathcal{U} L e_m , \psi_{n,\epsilon} \right>
\end{equation}
for $n \geq 1$, $\epsilon \in \{-1,+1\}$, and $1 \leq m \leq 2$. Introducing the following vectors and matrices:
\begin{equation*}
c_n(t) = \begin{bmatrix} c_{n,-1}(t) \\ c_{n,+1}(t) \end{bmatrix} ,
\quad
p_{d,n}(t) = \begin{bmatrix} p_{d,n,-1}(t) \\ p_{d,n,+1}(t) \end{bmatrix} ,
\end{equation*}
\begin{equation*}
\Lambda_n = \begin{bmatrix} \lambda_{n,-1} & 0 \\ 0 & \lambda_{n,+1} \end{bmatrix} ,
\quad
B_n = \begin{bmatrix} b_{n,-1,1} & b_{n,-1,2} \\ b_{n,+1,1} & b_{n,+1,2} \end{bmatrix} ,
\end{equation*}
and
\begin{equation}\label{eq: definition Mn}
M_n = \dfrac{\gamma}{\lambda_{n,-1}-\lambda_{n,+1}} \begin{bmatrix} 1 & \dfrac{k_{n,-1}}{k_{n,+1}} \\ -\dfrac{k_{n,+1}}{k_{n,-1}} & -1 \end{bmatrix} ,
\end{equation}
we obtain the following spectral decomposition:
\begin{align}
\dot{c}_n(t) & = \Lambda_n c_n(t) + M_n \{ c_n(t-h(t)) - c_n(t) \} \label{eq: spectral decomposition - prel bis} \\
& \phantom{=}\, + B_n \{u(t)+d_b(t)\} + p_{d,n}(t) \nonumber
\end{align}
for all $n \geq 1$ and $t \geq 0$. The initial condition is given by
\begin{equation}\label{eq: spectral decomposition IC - prel bis} 
c_n(\tau) = c_{\Phi,n}(\tau) \triangleq \begin{bmatrix} \left< \Phi(\tau) , \psi_{n,-1} \right> \\ \left< \Phi(\tau) , \psi_{n,+1} \right> \end{bmatrix}
\end{equation}
for $\tau\in[-h_M,0]$.

\begin{remark}
Even if the particularly selected lifting operator $L$ associated with $(\mathcal{A},\mathcal{B})$ is not unique, the coefficients $b_{n,\epsilon,m}$ given by (\ref{eq: coeff b_n}), and thus the resulting input matrices $B_n$, are actually independent of such a particular selection. See~\cite{lhachemi2019feedback} for a detailed explanation.
\end{remark}

\subsection{Finite-dimensional truncated model}\label{subsec: truncated model for control design}

We can now introduce a truncated model of the Euler-Bernoulli beam model (\ref{eq: beam model - PDE}-\ref{eq: beam model - IC2}) capturing a finite number of modes. For a given integer $N_0 \geq 1$, which will be discussed in the sequel, we introduce the following vectors and matrices:
\begin{subequations}
\begin{align}
Y(t) & = \begin{bmatrix} c_1(t)^\top & \ldots & c_{N_0}(t)^\top \end{bmatrix}^\top \in \mathbb{R}^{2 N_0} , \label{eq: truncated model - Y} \\
P_d(t) & = \begin{bmatrix} p_{d,1}(t)^\top & \ldots & p_{d,N_0}(t)^\top \end{bmatrix}^\top \in \mathbb{R}^{2 N_0} , \label{eq: truncated model - D} \\
Y_\Phi(\tau) & = \begin{bmatrix} c_{\Phi,1}(\tau)^\top & \ldots & c_{\Phi,N_0}(\tau)^\top \end{bmatrix}^\top \in \mathbb{R}^{2 N_0} , \label{eq: truncated model - Y,Phi} \\
A & = \mathrm{diag}(\Lambda_1,\ldots,\Lambda_{N_0}) \in \mathbb{R}^{2 N_0 \times 2 N_0} , \label{eq: truncated model - A} \\
B & = \begin{bmatrix} B_1^\top & \ldots & B_{N_0}^\top \end{bmatrix}^\top \in \mathbb{R}^{2 N_0 \times 2} , \label{eq: truncated model - B} \\
M & = \mathrm{diag}(M_1,\ldots,M_{N_0}) \in \mathbb{R}^{2 N_0 \times 2 N_0} . \label{eq: truncated model - M}
\end{align}
\end{subequations}
Then, we infer from (\ref{eq: spectral decomposition - prel bis}-\ref{eq: spectral decomposition IC - prel bis}) that the following ODE holds:
\begin{subequations}
\begin{align}
\dot{Y}(t) & = A Y(t) + M \{ Y(t-h(t)) - Y(t) \} \label{eq: spectral decomposition} \\
& \phantom{=}\; + B \{u(t)+d_b(t)\} + P_d(t) \nonumber \\
Y(\tau) & = Y_\Phi(\tau) , \quad \tau\in[-h_M,0] \label{eq: spectral decomposition IC}
\end{align}
\end{subequations}
for all $t \geq 0$.

Denoting by $B_{c1},B_{c2} \in\mathbb{R}^{2N_0}$ the first and second columns of the matrix $B$, respectively, we have the controllability property stated below. This result assesses the existence of a matrix $K \in \mathbb{R}^{2 \times 2 N_0}$ such that $A_\mathrm{cl} \triangleq A+BK$ is Hurwitz with desired pole placement. 

\begin{lemma}\label{lem: kalman condition}
For any given $N_0 \geq 1$, the pairs $(A,B)$, $(A,B_{c1})$, and $(A,B_{c2})$ satisfy the Kalman condition.
\end{lemma}

\textbf{Proof.}
Due to the diagonal nature of the matrix $A$ and the fact that the eigenvalues of $\mathcal{U}_0$ are simple, the PBH test~\cite{zhou1996robust} ensures the satisfaction of the Kalman condition provided $b_{n,\epsilon,m} \neq 0$. This is indeed the case because, based on (\ref{eq: coeff b_n}), straightforward computations show that $b_{n,\epsilon,1} = - n \pi C_{n,\epsilon} \neq 0$ and $b_{n,\epsilon,2} = (-1)^n n \pi C_{n,\epsilon} \neq 0$.
\qed

\begin{remark}\label{rem: number of control input}
The proposed approach captures the case of a single boundary control input (either at $x=0$ or at $x=1$). Indeed, one can impose $u_m = 0$ by setting the $m$-th line of the feedback gain $K$ as $0_{1 \times 2N_0}$. 
\end{remark}

\subsection{Control strategy and main result}

The proposed control strategy for stabilizing the Euler-Bernoulli beam model (\ref{eq: beam model - PDE}-\ref{eq: beam model - IC2}) consists in the two following steps. First, a feedback gain $K$ is computed such that the feedback law $u = KY$ exponentially stabilizes the truncated model (\ref{eq: spectral decomposition}-\ref{eq: spectral decomposition IC}) for some value $h_M > 0$ of the maximum amplitude of the admissible state-delay $h$. Specifically, for a given feedback gain $K$, we provide an LMI-based (sufficient) condition on the value of $h_M$. Second, we ensure that the design performed on the truncated model successfully stabilizes the original infinite-dimensional system (\ref{eq: beam model - PDE}-\ref{eq: beam model - IC2}) under the assumption that the number of modes $N_0$ of the truncated model is large enough. More precisely, introducing for any $\Phi \in \mathcal{C}^1([-h_M,0];\mathcal{H})$
\begin{equation*}
\Vert \Phi \Vert_{1,h_M} 
= \sqrt{ \Vert \Phi(0) \Vert^2 + \int_{-h_M}^{0} \Vert \dot{\Phi}(\tau) \Vert^2 \diff\tau} ,
\end{equation*}
the main result of this paper is stated as follows.

\begin{theorem}\label{thm: main theorem}
Let $N_0 \geq 1$ be such that $N_0 \geq \lfloor \beta^{1/4}/\pi \rfloor$ and\footnote{We recall that $\beta = \beta_0 + \gamma$ with $\alpha,\beta_0,\gamma$ the parameters of the studied Euler-Bernoulli beam model (\ref{eq: beam model - PDE}-\ref{eq: beam model - IC2}).}
\begin{equation}\label{eq: main thm - small gain condition}
\dfrac{60 \alpha^2 \gamma^2}{(\alpha^2 - 1) (N_0 + 1)^4 \pi^4 + \beta} \times \dfrac{1}{\lambda_{N_0+1,+1}^2} < 1 .
\end{equation}
Consider the matrices $A$ and $B$ defined by (\ref{eq: truncated model - A}-\ref{eq: truncated model - B}). Let $K \in \mathbb{R}^{2 \times 2 N_0}$ be such that $A_\mathrm{cl} = A + BK$ is Hurwitz and let $0 < h_m < h_M$ be such that
\begin{equation}\label{eq: main thm - LMI condition}
\begin{bmatrix}
A_\mathrm{cl}^\top P_2 + P_2^\top A_\mathrm{cl} & P_1 - P_2^\top + A_\mathrm{cl}^\top P_3 & h_M P_2^\top M \\
P_1 - P_2 + P_3^\top A_\mathrm{cl} & - P_3 - P_3^\top + h_M Q & h_M P_3^\top M \\
h_M M^\top P_2 & h_M M^\top P_3 & \; - h_M Q
\end{bmatrix} 
\prec 0 
\end{equation}
for some matrices $P_1,Q\in\mathbb{S}_{2N_0}^{+*}$, and $P_2, P_3 \in \mathbb{R}^{2N_0 \times 2N_0}$. Then, there exist constants $\kappa , C_0 , C_1 , C_2 > 0$ such that, for any distributed perturbation $d_d \in \mathcal{C}^0(\mathbb{R}_+;L^2(0,1))$, any boundary perturbation $d_b \in \mathcal{C}^1(\mathbb{R}_+;\mathbb{R}^2)$, any initial condition $\Phi \in \mathcal{C}^1([-h_M,0];\mathcal{H})$, and any delay $h \in \mathcal{C}^0(\mathbb{R}_+;\mathbb{R})$ with $h_m \leq h \leq h_M$, the mild solution $X\in\mathcal{C}^0(\mathbb{R}_+;\mathcal{H})$ of (\ref{eq: abstract form - DE}-\ref{eq: abstract form - IC}) with $u = K Y \in \mathcal{C}^1(\mathbb{R}_+;\mathbb{R}^2)$ satisfies, for all $t \geq 0$, 
\begin{align}
\Vert ( y(t,\cdot) , y_t(t,\cdot) ) \Vert 
& \leq C_0 e^{- \kappa t} \Vert (y_0,y_{t0}) \Vert_{1,h_M} \label{eq: main theorem - estimate y} \\
& \phantom{\leq}\; + C_1 \sup\limits_{\tau \in [0,t]} e^{- \kappa (t-\tau)}\Vert d_d(\tau) \Vert \nonumber \\
& \phantom{\leq}\; + C_2 \sup\limits_{\tau \in [0,t]} e^{- \kappa (t-\tau)}\Vert d_b(\tau) \Vert , \nonumber
\end{align}
with control input $\Vert u(t) \Vert \leq \frac{\Vert K \Vert}{\sqrt{m_R}} \Vert ( y(t,\cdot) , y_t(t,\cdot) ) \Vert$, where $( y(t,\cdot) , y_t(t,\cdot) ) = X(t)$ and $(y_0,y_{t0}) = \Phi$.
\end{theorem}

\begin{remark}
The condition (\ref{eq: main thm - small gain condition}) is always satisfied for $N_0$ large enough because its left hand side goes to zero as $N_0 \rightarrow + \infty$. Furthermore, because $A_\mathrm{cl}$ is Hurwitz, there always exists a $h_M > 0$ such that the LMI (\ref{eq: main thm - LMI condition}) is feasible. This latter result will be assessed in the sequel; see Lemma~\ref{lem: prel lemma feasability} for details. Consequently, conclusions of Theorem~\ref{thm: main theorem} always hold true for $N_0$ large enough and $h_M > 0$ small enough.
\end{remark}

\begin{remark}
From Remark~\ref{rem: properties eigenvalues}, the condition $N_0 \geq \lfloor \beta^{1/4}/\pi \rfloor$ ensures that all the unstable modes of the Euler-Bernoulli beam model\footnote{More specifically all the unstable eigenvalues of the disturbance-free operator $\mathcal{U}_0$.} (\ref{eq: beam model - PDE}-\ref{eq: beam model - IC2}) are captured by the truncated model (\ref{eq: spectral decomposition}-\ref{eq: spectral decomposition IC}). In particular, $\lambda_{n,\epsilon} \leq \lambda_{N_0+1,+1} < 0$ for all $n \geq N_0 +1$ and $\epsilon\in\{-1,+1\}$.
\end{remark}

\begin{remark}\label{rem: ISS-like estimate}
Note that the estimate (\ref{eq: main theorem - estimate y}) is not an ISS estimate in strict form because the system trajectory and the initial condition are not evaluated with the same norm. This issue can be overcome as follows. Assume that the assumptions of Theorem~\ref{thm: main theorem} hold but with the LMI condition (\ref{eq: main thm - LMI condition}) replaced by the small gain condition\footnote{By a continuity argument in $h_M = 0$, there always exists $h_M > 0$ such that the small gain condition (\ref{eq: main thm - ISS in strict form}) holds.}
\begin{equation}\label{eq: main thm - ISS in strict form}
C_\lambda \Vert M \Vert \left( e^{\Vert A_\mathrm{cl} \Vert h_M} - e^{- \lambda h_M} \right) < \lambda ,
\end{equation}
where $\lambda > 0$ and $C_\lambda \geq 1$ are any constants such that $\Vert e^{A_\mathrm{cl} t} \Vert \leq C_\lambda e^{-\lambda t}$ for all $t \geq 0$. Then there exist constants $\overline{\kappa},\overline{C}_0,\overline{C}_1,\overline{C}_2 > 0$, independent of $\Phi,d_d,d_b,h$, such that the mild solutions satisfy, for all $t \geq 0$, 
\begin{align}
\Vert ( y(t,\cdot) , y_t(t,\cdot) ) \Vert 
& \leq \overline{C}_0 e^{- \overline{\kappa} t} \sup\limits_{\tau\in[-h_M,0]} \Vert (y_0(\tau,\cdot),y_{t0}(\tau,\cdot)) \Vert \label{eq: main theorem - ISS in strict form - estimate y} \\
& \phantom{\leq}\; + \overline{C}_1 \sup\limits_{\tau \in [0,t]} e^{- \overline{\kappa} (t-\tau)}\Vert d_d(\tau) \Vert \nonumber \\
& \phantom{\leq}\; + \overline{C}_2 \sup\limits_{\tau \in [0,t]} e^{- \overline{\kappa} (t-\tau)}\Vert d_b(\tau) \Vert , \nonumber
\end{align}
with control input $\Vert u(t) \Vert \leq \frac{\Vert K \Vert}{\sqrt{m_R}} \Vert ( y(t,\cdot) , y_t(t,\cdot) ) \Vert$ and where, comparing to Theorem~\ref{thm: main theorem}, the regularity assumption on the initial condition is weakened to $\Phi\in\mathcal{C}^0(\mathbb{R}_+;\mathcal{H})$.

The ISS estimate provided by (\ref{eq: main theorem - ISS in strict form - estimate y}) is sharper than the ISS-like estimate (\ref{eq: main theorem - estimate y}). However, the small gain condition (\ref{eq: main thm - ISS in strict form}) provides, in general, numerical values of $h_M$ that are significantly lower than the ones provided by the LMI condition (\ref{eq: main thm - LMI condition}). In particular, any $h_M > 0$ satisfying the small gain condition (\ref{eq: main thm - ISS in strict form}) is such that
\begin{equation}\label{eq: main thm - ISS in strict form - upperbound}
h_M < \dfrac{1}{\Vert A_\mathrm{cl} \Vert} \log\left( 1 + \dfrac{\mu_M(A_\mathrm{cl})}{\Vert M \Vert} \right) ,
\end{equation}
where $0 < \mu_M(A_\mathrm{cl}) = - \max\left\{ \operatorname{Re}(\lambda) \,:\, \lambda \in \mathrm{sp}_\mathbb{C}(A_\mathrm{cl}) \right\}$.
\end{remark}

\begin{remark}
From the proof of Theorem~\ref{thm: main theorem} reported in the sequel (see in particular Lemma~\ref{lem: infinite-dim part negelected in the design}), we can also derive the following result regarding the behavior of the open-loop system (i.e. $u=0$). Assume that $\beta < \pi^4$ and that the small gain condition (\ref{eq: main thm - small gain condition}) holds for $N_0 = 0$. The former condition ensures that all the eigenvalues of the disturbance-free operator $\mathcal{U}_0$ are stable. Recalling that $\beta = \beta_0+\gamma$, these two conditions are satisfied for small enough values of $\beta_0 \geq 0$ and $\gamma>0$. Let $0 < h_m < h_M$ be arbitrary. Then there exist constants $\kappa , C_0 , C_1 , C_2 > 0$ such that, for any $d_d,d_b,\Phi,h$ satisfying the regularity assumptions of Theorem~\ref{thm: main theorem}, the mild solution of (\ref{eq: abstract form - DE}-\ref{eq: abstract form - IC}) with $u = 0$ satisfies the estimate (\ref{eq: main theorem - estimate y}).
\end{remark}

From Remark~\ref{rem: Pi bounded}, we have for any $y = (y_1,y_2) \in \mathcal{H}$ that $\vert y_1(x) \vert \leq \Vert y \Vert$ for all $x \in [0,1]$. This yields the following result regarding the uniform convergence of the displacements of the beam.

\begin{corollary}
Assume that the assumptions of Theorem~\ref{thm: main theorem} hold true. Then, we have for all $t \geq 0$
\begin{align}
\sup\limits_{x \in [0,1]} \vert y(t,x) \vert 
& \leq C_0 e^{- \kappa t} \Vert (y_0(t,\cdot),y_{t0}(t,\cdot)) \Vert_{1,h_M} \label{eq: cor main theorem - estimate displacement sup norm}\\
& \phantom{\leq}\; + C_1 \sup\limits_{\tau \in [0,t]} e^{- \kappa (t-\tau)} \Vert d_d(\tau) \Vert \nonumber \\
& \phantom{\leq}\; + C_2 \sup\limits_{\tau \in [0,t]} e^{- \kappa (t-\tau)}\Vert d_b(\tau) \Vert . \nonumber
\end{align}
\end{corollary}

\subsection{Well-posedness of the closed-loop system}\label{subsec: well-posedness}

The well-posedness of the closed-loop system (\ref{eq: abstract form - DE}-\ref{eq: abstract form - IC}) with $u = K Y$ in terms of mild solutions is assessed by the following lemma whose proof is placed in Appendix~\ref{annex: well-posedness}. 

\begin{lemma}\label{lem: existence mild solutions for the closed-loop system}
Let $0 < h_m < h_M$, $N_0 \geq 1$, and $K \in \mathbb{R}^{2 \times 2 N_0}$ be arbitrary. Let $d_d \in \mathcal{C}^0(\mathbb{R}_+;L^2(0,1))$, $d_b \in \mathcal{C}^1(\mathbb{R}_+;\mathbb{R}^2)$, $\Phi \in \mathcal{C}^0([-h_M,0];\mathcal{H})$, and $h \in \mathcal{C}^0(\mathbb{R}_+;\mathbb{R})$ with $h_m \leq h \leq h_M$ be arbitrary. Then, there exists a unique mild solution $X\in\mathcal{C}^0(\mathbb{R}_+;\mathcal{H})$ of (\ref{eq: abstract form - DE}-\ref{eq: abstract form - IC}) with control input $u = K Y \in \mathcal{C}^1(\mathbb{R}_+;\mathbb{R}^2)$.
\end{lemma}

The above result establishes the validity of the spectral decomposition and associated truncated model, as reported in Subsections~\ref{subsec: spectral decomposition} and~\ref{subsec: truncated model for control design}, for the abstract boundary control system (\ref{eq: abstract form - DE}-\ref{eq: abstract form - IC}) when placed in closed loop with $u = K Y$.

\begin{remark}
It is possible to extend the result stated in Theorem~\ref{thm: main theorem} to any boundary disturbance with relaxed regularity assumption $d_b \in \mathcal{C}^0(\mathbb{R}_+;\mathbb{R}^2)$ by considering the concept of weak solution as introduced in~\cite{lhachemi2018iss}. Specifically, for any $d_d \in \mathcal{C}^0(\mathbb{R}_+;L^2(0,1))$, $d_b \in \mathcal{C}^0(\mathbb{R}_+;\mathbb{R}^2)$, $\Phi \in \mathcal{C}^0([-h_M,0];\mathcal{H})$, and $h \in \mathcal{C}^0(\mathbb{R}_+;\mathbb{R})$ with $0 < h_m \leq h \leq h_M$, we say that $X\in\mathcal{C}^0(\mathbb{R}_+;\mathcal{H})$ is a weak solution of the closed-loop system (\ref{eq: abstract form - DE}-\ref{eq: abstract form - IC}) with $u = KY$ if for any $T > 0$ and any $z \in \mathcal{C}^0([0,T];D(\mathcal{A}_0^*))\cap\mathcal{C}^1([0,T];\mathcal{H})$ such that $\mathcal{A}_0^* z \in\mathcal{C}^0([0,T];\mathcal{H})$ and $z(T)=0$,
\begin{align}
& \int_0^T \left< X(t) , \mathcal{A}_0^* z(t) + \dfrac{\mathrm{d}z}{\mathrm{d}t}(t) \right> \diff t \nonumber \\
& = - \left< \Phi(0) , z(0) \right>
+ \int_0^T \left< L w(t) , \mathcal{A}_0^* z(t) \right> \diff t \label{eq: def weak solution} \\
& \phantom{=}\; - \int_0^T \left< \mathcal{A}Lw(t) + \gamma\Pi X(t-h(t)) + p_d(t) , z(t) \right> \diff t , \nonumber
\end{align}
where $w=KY+d_b$ with $Y$ defined by (\ref{eq: truncated model - Y}) and with the initial condition $X(\tau) = \Phi(\tau)$ for all $\tau\in[-h_M,0]$. First, it can be shown that any mild solution is a weak solution. Second, using the test functions $z(t) = (t-T)\psi_{n,\epsilon}$ over $[0,T]$ and the identity $\mathcal{U}_0 = \mathcal{A}_0 + \gamma\Pi$ with $\Pi$ bounded, it can be shown the uniqueness of the weak solutions. Finally, under the assumptions of Theorem~\ref{thm: main theorem}, as the mild solutions of the closed-loop system (\ref{eq: abstract form - DE}-\ref{eq: abstract form - IC}) with $u = KY$ satisfy the estimate (\ref{eq: main theorem - estimate y}), a density argument similar to the one reported in~\cite[Proof of Thm.~3]{lhachemi2018iss} can be used to prove the existence of the weak solutions associated with any\footnote{In the context of Remark~\ref{rem: ISS-like estimate}, these conclusions also hold, when substituting the ISS-like estimate (\ref{eq: main theorem - estimate y}) by the ISS estimate (\ref{eq: main theorem - ISS in strict form - estimate y}), for all initial condition $\Phi \in \mathcal{C}^0([-h_M,0])$.} $d_d \in \mathcal{C}^0(\mathbb{R}_+;L^2(0,1))$, $d_b \in \mathcal{C}^0(\mathbb{R}_+;\mathbb{R}^2)$, $\Phi \in \mathcal{C}^1([-h_M,0])$, and $h \in \mathcal{C}^0(\mathbb{R}_+;\mathbb{R})$ such that $h_m \leq h \leq h_M$. Moreover, these weak solutions satisfy the ISS-like estimate (\ref{eq: main theorem - estimate y}).
\end{remark}

\section{Stability of the closed-loop finite-dimensional truncated model}\label{sec: stab CL finite-dim}

The objective of this section is to derive values of $h_M > 0$ that ensures the existence of an exponential ISS-like estimate for the truncated model (\ref{eq: spectral decomposition}-\ref{eq: spectral decomposition IC}) with $u=KY$ where $K \in \mathbb{R}^{2 \times 2 N_0}$ is selected such that $A_\mathrm{cl} = A+BK$ is Hurwitz.

\subsection{Preliminary lemmas}

For $h_M > 0$, following~\cite[Chap.~4, Sec.~1.3]{kolmanovskii2012applied}, we introduce $W$ the space of absolutely continuous functions $f : [-h_M,0] \rightarrow \mathbb{R}^n$ with square-integrable derivative endowed with the norm $\Vert f \Vert_{W} \triangleq \sqrt{ \Vert f(0) \Vert^2 + \int_{-h_M}^0 \Vert \dot{f}(\tau) \Vert^2 \diff \tau}$. The following Lemma is an exponential ISS-like version of a disturbance-free asymptotic stability result reported in~\cite{fridman2006new,fridman2014tutorial}.

\begin{lemma}\label{lem: prel lemma}
Let $F,G \in \mathbb{R}^{n \times n}$ and $0 < h_m < h_M$ be given. Assume that there exist $\kappa > 0$, $P_1,Q\in\mathbb{S}_n^{+*}$, and $P_2, P_3 \in \mathbb{R}^{n \times n}$ such that $\Theta(h_M,\kappa) \prec 0$ with
{\small
\begin{align}
& \Theta(h_M,\kappa) = \label{eq: prel lemma - LMI condition} \\
& \begin{bmatrix}
2 \kappa P_1 + F^\top P_2 + P_2^\top F & P_1 - P_2^\top + F^\top P_3 & h_M P_2^\top G \\
P_1 - P_2 + P_3^\top F & - P_3 - P_3^\top + h_M Q & h_M P_3^\top G \\
h_M G^\top P_2 & h_M G^\top P_3 & \; - h_M e^{-2 \kappa h_M} Q
\end{bmatrix} . \nonumber 
\end{align}
}
Then, there exist constants $\gamma_0,\gamma_1>0$ such that, for any  $h \in \mathcal{C}^0(\mathbb{R}_+;\mathbb{R}_+)$ with $h_m \leq h \leq h_M$, and any $d \in \mathcal{C}^0(\mathbb{R}_+;\mathbb{R}^n)$, the trajectory $x$ of
\begin{align*}
\dot{x}(t) & = F x(t) + G \left\{ x(t-h(t)) - x(t) \right\} + d(t) \\
x(\tau) & = x_0(\tau) , \quad \tau \in [-h_M,0] 
\end{align*}
for $t \geq 0$ and with initial condition $x_0 \in W$ satisfies 
\begin{equation*}
\Vert x(t) \Vert \leq \gamma_0 e^{- \kappa t} \Vert x_0 \Vert_W + \gamma_1 \sup\limits_{\tau\in[0,t]} e^{-\kappa(t-\tau)} \Vert d(\tau) \Vert 
\end{equation*}
for all $t \geq 0$. 
\end{lemma}

\textbf{Proof.} 
Let $\kappa > 0$, $P_1,Q\in\mathbb{S}_n^{+*}$, and $P_2, P_3 \in \mathbb{R}^{n \times n}$ be such that $\Theta(h_M,\kappa) \prec 0$. By a continuity argument, we select $\sigma > \kappa$ such that $\Theta(h_M,\sigma) \prec 0$. We consider the Lyapunov-Krasovskii functional $V(t) = V_1(t) + V_2(t)$ with 
\begin{align*}
V_1(t) & = x(t)^\top P_1 x(t) , \\
V_2(t) & = \int_{-h_M}^{0} \int_{t+\theta}^{t} e^{2\sigma(s-t)} \dot{x}(s)^\top Q \dot{x}(s) \diff s \diff \theta .
\end{align*}
Then, the computation of the time derivative yields 
\begin{align}
\dot{V}(t)
& = 2 x(t)^\top P_1 \dot{x}(t) + h_M \dot{x}(t)^\top Q \dot{x}(t) - 2 \sigma V_2(t) \label{eq: computation dot_V} \\
& \phantom{=}\, - \int_{-h_M}^{0} e^{2 \sigma \theta} \dot{x}(t+\theta)^\top Q \dot{x}(t+\theta) \diff \theta , \nonumber
\end{align}
for all $t \geq 0$. Now, noting that 
\begin{equation}
\dot{x}(t) = F x(t) - G \int_{t-h(t)}^{t} \dot{x}(\tau) \diff \tau + d(t) \label{eq: dynamics in function of dot_x}
\end{equation}
for all $t \geq 0$,  and introducing
\begin{equation*}
P = \begin{bmatrix} P_1 & 0 \\ P_2 & P_3 \end{bmatrix} 
\end{equation*} 
where $P_2, P_3 \in \mathbb{R}^{n \times n}$ are ``slack variables''~\cite{fridman2006new,fridman2014tutorial}, we obtain that
\begin{align}
& x(t)^\top P_1 \dot{x}(t) \nonumber \\
& \overset{(\ref{eq: dynamics in function of dot_x})}{=}  
\begin{bmatrix}
x(t) \\ \dot{x}(t)
\end{bmatrix}^\top
P^\top
\begin{bmatrix}
\dot{x}(t) \\ -\dot{x}(t) + F x(t) - G \int_{t-h(t)}^{t} \dot{x}(\tau) \diff \tau + d(t)
\end{bmatrix} \nonumber  \\
& = 
\begin{bmatrix}
x(t) \\ \dot{x}(t)
\end{bmatrix}^\top
P^\top
\begin{bmatrix}
0 & I \\ F & -I 
\end{bmatrix}
\begin{bmatrix}
x(t) \\ \dot{x}(t)
\end{bmatrix}  
+ \begin{bmatrix}
x(t) \\ \dot{x}(t)
\end{bmatrix}^\top
P^\top
\begin{bmatrix}
0 \\ d(t)
\end{bmatrix} \label{eq: coumputation x'*P1*dot_x} \\
& \phantom{=}\; + \int_{t-h(t)}^{t}  
\begin{bmatrix}
x(t) \\ \dot{x}(t)
\end{bmatrix}^\top
P^\top
\begin{bmatrix}
0 \\ - G
\end{bmatrix}
\dot{x}(\tau) \diff \tau . \nonumber 
\end{align}
We estimate the two last terms. First, as $2 a^\top b \leq \Vert a \Vert^2 + \Vert b \Vert^2$ for any $a,b \in \mathbb{R}^n$, we obtain that
\begin{align*}
& 2 \begin{bmatrix}
x(t) \\ \dot{x}(t)
\end{bmatrix}^\top
P^\top
\begin{bmatrix}
0 \\ - G
\end{bmatrix}
\dot{x}(\tau) \\
& \quad = 2 \left( e^{- \sigma (\tau-t)} Q^{-1/2} 
\begin{bmatrix}
0 \\ - G
\end{bmatrix}^\top
P
\begin{bmatrix}
x(t) \\ \dot{x}(t)
\end{bmatrix} \right)^\top \\
& \quad \phantom{=}\; \times
\left( e^{\sigma (\tau-t)} Q^{1/2} \dot{x}(\tau) \right) \\
& \quad \leq 
e^{-2 \sigma (\tau-t)}
\begin{bmatrix}
x(t) \\ \dot{x}(t)
\end{bmatrix}^\top
P^\top
\begin{bmatrix}
0 \\ G
\end{bmatrix}
Q^{-1}
\begin{bmatrix}
0 \\ G
\end{bmatrix}^\top
P
\begin{bmatrix}
x(t) \\ \dot{x}(t)
\end{bmatrix} \\ 
& \quad \phantom{\leq}\,  + e^{2 \sigma (\tau-t)} \dot{x}(\tau)^\top Q \dot{x}(\tau) .
\end{align*}
Noting that 
\begin{equation*}
0 \leq \int_{t-h(t)}^t e^{-2 \sigma \tau} \leq h(t) e^{- 2 \sigma (t-h(t))} \leq h_M e^{2 \sigma h_M} e^{-2\sigma t} , 
\end{equation*}
we infer that
\begin{align*}
& 2 \int_{t-h(t)}^{t}  
\begin{bmatrix}
x(t) \\ \dot{x}(t)
\end{bmatrix}^\top
P^\top
\begin{bmatrix}
0 \\ - G
\end{bmatrix}
\dot{x}(\tau) \diff \tau \\
& \;\;\leq h_M e^{2 \sigma h_M} \begin{bmatrix}
x(t) \\ \dot{x}(t)
\end{bmatrix}^\top
P^\top
\begin{bmatrix}
0 \\ G
\end{bmatrix}
Q^{-1}
\begin{bmatrix}
0 \\ G
\end{bmatrix}^\top
P
\begin{bmatrix}
x(t) \\ \dot{x}(t)
\end{bmatrix} \\
& \phantom{\;\;\leq}\; + \int_{t-h(t)}^{t} e^{2 \sigma (\tau-t)} \dot{x}(\tau)^\top Q \dot{x}(\tau) \diff \tau .
\end{align*}
We now estimate the term related to the disturbance input $d(t)$. For any constant $\epsilon > 0$, which will be specified later, the use of Young's inequality yields
\begin{align*}
& \begin{bmatrix}
x(t) \\ \dot{x}(t)
\end{bmatrix}^\top
P^\top
\begin{bmatrix}
0 \\ d(t)
\end{bmatrix} \\
& \quad = x(t)^\top P_2^\top d(t) + \dot{x}(t)^\top P_3^\top d(t) \\
& \quad \leq \dfrac{\epsilon}{2} \left( \Vert x(t) \Vert^2 + \Vert \dot{x}(t) \Vert^2 \right) 
+ \dfrac{\Vert P_2^\top \Vert^2 + \Vert P_3^\top \Vert^2}{2\epsilon} \Vert d(t) \Vert^2 .
\end{align*}
Combining the two latter estimates with (\ref{eq: computation dot_V}) and (\ref{eq: coumputation x'*P1*dot_x}), we deduce that
\begin{align*}
& \dot{V}(t) + 2 \sigma V(t)\\
& \quad \leq 2 \sigma V_1(t) +
2
\begin{bmatrix}
x(t) \\ \dot{x}(t)
\end{bmatrix}^\top
P^\top
\begin{bmatrix}
0 & I \\ F & -I 
\end{bmatrix}
\begin{bmatrix}
x(t) \\ \dot{x}(t)
\end{bmatrix} \\
& \quad \phantom{\leq}\, + h_M e^{2 \sigma h_M}
\begin{bmatrix}
x(t) \\ \dot{x}(t)
\end{bmatrix}^\top
P^\top
\begin{bmatrix}
0 \\ G
\end{bmatrix}
Q^{-1}
\begin{bmatrix}
0 \\ G
\end{bmatrix}^\top
P
\begin{bmatrix}
x(t) \\ \dot{x}(t)
\end{bmatrix} \\
& \quad \phantom{\leq}\, + h_M \dot{x}(t)^\top Q \dot{x}(t)
- \int_{t-h_M}^{t-h(t)} e^{2 \sigma (\tau-t)} \dot{x}(\tau)^\top Q \dot{x}(\tau) \diff \tau \\
& \quad \phantom{\leq}\, + \epsilon \begin{bmatrix}
x(t) \\ \dot{x}(t)
\end{bmatrix}^\top \begin{bmatrix}
x(t) \\ \dot{x}(t)
\end{bmatrix}
+ \dfrac{\Vert P_2^\top \Vert^2 + \Vert P_3^\top \Vert^2}{\epsilon} \Vert d(t) \Vert^2 \\
& \quad \leq 
\begin{bmatrix}
x(t) \\ \dot{x}(t)
\end{bmatrix}^\top
\left\{
\Xi
+ \epsilon I
\right\}
\begin{bmatrix}
x(t) \\ \dot{x}(t)
\end{bmatrix} 
+ \dfrac{\Vert P_2^\top \Vert^2 + \Vert P_3^\top \Vert^2}{\epsilon} \Vert d(t) \Vert^2 ,
\end{align*}
where it has been used the fact that the integral term is always non negative because the integrand is non negative and $h(t) \leq h_M$, and where
\begin{align*}
\Xi & \triangleq
P^\top
\begin{bmatrix}
0 & I \\ F & -I
\end{bmatrix}
+
\begin{bmatrix}
0 & I \\ F & -I
\end{bmatrix}^\top
P
+ 
\begin{bmatrix}
2 \sigma P_1 & 0 \\ 0 & h_M Q 
\end{bmatrix} \\
& \phantom{=}\; + h_M e^{2 \sigma h_M}
P^\top
\begin{bmatrix}
0 \\ G
\end{bmatrix}
Q^{-1}
\begin{bmatrix}
0 \\ G
\end{bmatrix}^\top
P .
\end{align*}
From $\Theta(h_M,\sigma) \prec 0$, the Schur complement shows $\Xi \prec 0$. Thus, we can select $\epsilon > 0$ small enough, independently of $h$ and $x_0$, such that $\Xi + \epsilon I \preceq 0$. This yields, for all $t \geq 0$,
\begin{equation*}
\dot{V}(t) + 2\sigma V(t) \leq \delta \Vert d(t) \Vert^2 
\end{equation*}
with $\delta \triangleq (\Vert P_2^\top \Vert^2 + \Vert P_3^\top \Vert^2)/\epsilon$, giving
\begin{equation*}
V(t) \leq e^{-2\sigma t} V(0) + \delta e^{-2\sigma t} \int_0^t e^{2\sigma\tau} \Vert d(\tau) \Vert^2 \diff\tau .
\end{equation*}
Introducing $\eta = \kappa/\sigma \in(0,1)$, successive estimations show
\begin{align*}
& e^{-2\sigma t} \int_0^t e^{2\sigma\tau} \Vert d(\tau) \Vert^2 \diff\tau \\
& \quad = e^{-2\sigma t} \int_0^t e^{2\sigma(1-\eta)\tau} \times e^{2\sigma\eta\tau} \Vert d(\tau) \Vert^2 \diff\tau \\
& \quad \leq e^{-2\sigma t} \dfrac{e^{2\sigma(1-\eta)t}}{2\sigma(1-\eta)} \sup\limits_{\tau\in[0,t]} e^{2\sigma\eta\tau} \Vert d(\tau) \Vert^2 \\
& \quad \leq \dfrac{1}{2(\sigma-\kappa)} \sup\limits_{\tau\in[0,t]} e^{-2\kappa(t-\tau)} \Vert d(\tau) \Vert^2 .
\end{align*}
Hence, using again $\kappa < \sigma$, we obtain the following estimate
\begin{equation*}
V(t) \leq e^{-2\kappa t} V(0) + \dfrac{\delta}{2(\sigma-\kappa)} \sup\limits_{\tau\in[0,t]} e^{-2\kappa(t-\tau)} \Vert d(\tau) \Vert^2 
\end{equation*}
for all $t \geq 0$. Now the claimed conclusion follows from the fact that $\lambda_\mathrm{m}(P_1) \Vert x(t) \Vert^2 \leq  V(t) \leq \max\left( \lambda_\mathrm{M}(P_1) , h_M \lambda_\mathrm{M}(Q) \right) \Vert x(t+\cdot) \Vert_W^2$ for all $t \geq 0$ . \qed

Assuming that the matrix $F$ is Hurwitz, the following lemma ensures the feasability of the LMI $\Theta(h_M,\kappa) \prec 0$ for sufficiently small values of $h_M,\kappa > 0$.

\begin{lemma}\label{lem: prel lemma feasability}
Let $M,N \in \mathbb{R}^{n \times n}$ with $M$ Hurwitz be given. Then there exist $h_M,\kappa > 0$ such that $\Theta(h_M,0) \prec 0$ and $\Theta(h_M,\kappa) \prec 0$.
\end{lemma}

\textbf{Proof.}
By a continuity argument, it is sufficient to show the existence of $h_M > 0$ such that $\Theta(h_M,0) \prec 0$. To do so, let $P_2\in\mathbb{S}_n^{+*}$ be the unique solution of $F^\top P_2 + P_2 F = - I$. We define the matrices $P_1=2P_2\in\mathbb{S}_n^{+*}$, $P_3 = -(F^{-1})^\top P_2$, and $Q = I\in\mathbb{S}_n^{+*}$. Then we have $\Theta(h_M,0) \prec 0$ if and only if
\begin{equation*}
\begin{bmatrix}
-I & 0 & \; h_M P_2 G \\
0 & - S_3 + h_M I & \; h_M P_3^\top G \\
h_M G^\top P_2 & h_M G^\top P_3 & \; - h_M I
\end{bmatrix} 
\prec 0 
\end{equation*}
with $S_3 = P_3+P_3^\top = (F^{-1})^\top F^{-1} \succ 0$. Noting that $- h_M I \prec 0$, the Schur complement shows that the above LMI holds if and only if
\begin{equation}\label{eq: coro existence hM>0 - LMI condition}
\begin{bmatrix}
-I & 0 \\
0 & - S_3 + h_M I \\
\end{bmatrix}
+ h_M
\begin{bmatrix}
P_2 G \\
P_3^\top G 
\end{bmatrix} 
\begin{bmatrix}
P_2 G \\
P_3^\top G 
\end{bmatrix}^\top
\prec 0 .
\end{equation}
As $- \mathrm{diag}(I,S_3) \prec 0$, a continuity argument in $h_M = 0$ shows that (\ref{eq: coro existence hM>0 - LMI condition}) holds for sufficiently small $h_M > 0$. This completes the proof.
\qed

\subsection{Application to the truncated model}

We can now apply the results of the previous subsection to the truncated model (\ref{eq: spectral decomposition}-\ref{eq: spectral decomposition IC}) of the studied Euler-Bernoulli beam. In particular, in view of Lemmas~\ref{lem: kalman condition}, \ref{lem: prel lemma}, and~\ref{lem: prel lemma feasability}, and using the estimate
\begin{align*}
& \sup\limits_{\tau\in[0,t]} e^{-\kappa(t-\tau)} \Vert B d_b(\tau) + P_d(\tau) \Vert \\
& \leq
\Vert B \Vert \sup\limits_{\tau\in[0,t]} e^{-\kappa(t-\tau)} \Vert d_b(\tau) \Vert
+ \sup\limits_{\tau\in[0,t]} e^{-\kappa(t-\tau)} \Vert P_d(\tau) \Vert  ,
\end{align*}
we obtain the following result.

\begin{lemma}\label{lem: stab truncated model}
Let $N_0 \geq 1$ be given. Consider the matrices $A$ and $B$ defined by (\ref{eq: truncated model - A}-\ref{eq: truncated model - B}). Let $K \in \mathbb{R}^{2 \times {2 N_0}}$ be such that $A_\mathrm{cl} = A+BK$ is Hurwitz. Let $0<h_m<h_M$ and $\sigma>0$ be such that $\Theta(h_M,\sigma) \prec 0$ where $\Theta(h_M,\sigma)$ is defined by (\ref{eq: prel lemma - LMI condition}) with $F = A_\mathrm{cl}$ and $G = M$ given by (\ref{eq: truncated model - M}). Then, there exist constants $C_3, C_4 , C_5 > 0$ such that, for all $Y_\Phi \in \mathcal{C}^1([-h_M,0];\mathbb{R}^{N_0})$, $h \in \mathcal{C}^0(\mathbb{R}_+;\mathbb{R})$ with $h_m \leq h \leq h_M$, $P_d \in \mathcal{C}^0(\mathbb{R}_+;\mathbb{R}^{N_0})$, and $d_b \in \mathcal{C}^0(\mathbb{R}_+;\mathbb{R}^2)$, the trajectory $Y$ of (\ref{eq: spectral decomposition}-\ref{eq: spectral decomposition IC}) with command input $u = KY$ satisfies 
\begin{align}
\Vert Y(t) \Vert 
& \leq C_3 e^{- \sigma t} \Vert Y_\Phi \Vert_W + C_4 \sup\limits_{\tau \in [0,t]} e^{- \sigma(t-\tau)} \Vert P_d(\tau) \Vert \nonumber \\
& \phantom{\leq}\; C_5 \sup\limits_{\tau \in [0,t]} e^{- \sigma(t-\tau)} \Vert d_b(\tau) \Vert \label{eq: exp stab finite-dim part}
\end{align}
for all $t \geq 0$.
\end{lemma}

\begin{remark}
As discussed in Remark~\ref{rem: ISS-like estimate}, it is possible to transform the ISS-like estimate (\ref{eq: exp stab finite-dim part}) of Lemma~\ref{lem: stab truncated model} into an ISS estimate in strict form by replacing the LMI condition (\ref{eq: prel lemma - LMI condition}) with the small gain condition (\ref{eq: main thm - ISS in strict form}). This result relies on a small gain version of Lemma~\ref{lem: prel lemma} which can be derived similarly to~\cite[Thm.~2.5]{karafyllis2013delay}.
\end{remark}

\section{Exponential input-to-state stability of the closed-loop infinite-dimensional system}\label{sec: stab CL infinite-dim}

Now we complete the stability analysis of the infinite-dimensional Euler Bernoulli beam model (\ref{eq: beam model - PDE}-\ref{eq: beam model - IC2}) when placed in closed loop with the feedback law $u = KY$.

\subsection{Stability of the infinite-dimensional part of the system neglected in the control design}

As the control design has been performed on the truncated model (\ref{eq: spectral decomposition}-\ref{eq: spectral decomposition IC}) that only captures a finite number of modes of the original infinite-dimensional system (\ref{eq: beam model - PDE}-\ref{eq: beam model - IC2}), one has to ensure that the stability of the residual infinite-dimensional system is preserved.

\begin{lemma}\label{lem: infinite-dim part negelected in the design}
Let $0 < h_m < h_M$ and $\sigma,C_6,C_7,C_8>0$ be arbitrary. Let $N_0 \geq \lfloor \beta^{1/4}/\pi \rfloor$ be such that
\begin{equation}\label{eq: thm infinite dim residual - small gain condition}
\dfrac{60 \alpha^2 \gamma^2}{(\alpha^2 - 1) (N_0 + 1)^4 \pi^4 + \beta} \times \dfrac{1}{\lambda_{N_0+1,+1}^2} < 1 .
\end{equation}
Then, there exist constants $\kappa\in(0,\sigma)$ and $C_9 , C_{10},C_{11} > 0$ such that, for any $\Phi \in \mathcal{C}^1([-h_M,0];\mathcal{H})$, $d_d \in \mathcal{C}^0(\mathbb{R}_+;L^2(0,1))$, $d_b \in \mathcal{C}^1(\mathbb{R}_+;\mathbb{R}^2)$, $h \in \mathcal{C}^0(\mathbb{R}_+;\mathbb{R})$ such that $h_m \leq h \leq h_M$, and $u \in \mathcal{C}^1(\mathbb{R}_+;\mathbb{R}^2)$ such that, for all $t \geq 0$,
\begin{align}
\Vert u(t) \Vert 
& \leq C_6 e^{-\sigma t} \Vert \Phi \Vert_{1,h_M} + C_7 \sup\limits_{\tau \in [0,t]} e^{- \sigma(t-\tau)} \Vert d_d(\tau) \Vert \nonumber \\
& \phantom{\leq}\; + C_8 \sup\limits_{\tau \in [0,t]} e^{- \sigma(t-\tau)} \Vert d_b(\tau) \Vert , \label{eq: assumed estimate on u}
\end{align}
the mild solution $X\in\mathcal{C}^0(\mathbb{R}_+;\mathcal{H})$ of (\ref{eq: abstract form - DE}-\ref{eq: abstract form - IC}) satisfies
\begin{align}
\sum\limits_{\substack{n \geq N_0 + 1 \\ \epsilon\in\{-1,+1\}}} \vert c_{n,\epsilon}(t) \vert^2 
& \leq C_9 e^{- 2 \kappa t} \Vert \Phi \Vert_{1,h_M}^2 \label{eq: exp stab infinite-dim part} \\
& \phantom{\leq}\; + C_{10} \sup\limits_{\tau \in [0,t]} e^{- 2 \kappa (t-\tau)} \Vert d_d(\tau) \Vert^2 \nonumber \\
& \phantom{\leq}\; + C_{11} \sup\limits_{\tau \in [0,t]} e^{- 2 \kappa (t-\tau)} \Vert d_b(\tau) \Vert^2 \nonumber 
\end{align}
for all $t \geq 0$, where $c_{n,\epsilon}(t) = \left< X(t) , \psi_{n,\epsilon} \right>$.
\end{lemma}

The proposed proof is inspired by the small gain analysis reported in~\cite{karafyllis2013delay} in the context of finite-dimensional systems.

\textbf{Proof.}
We define $\eta = - \lambda_{N_{0}+1,+1}/2 > 0$. Throughout this proof, we always consider integers $n$ such that $n \geq N_0 + 1$, which ensures that $\lambda_{n,\epsilon} \leq \lambda_{N_0 +1,+1} = - 2 \eta < 0$ for all $\epsilon\in\{-1,+1\}$. Let $\kappa \in (0,\min(\eta,\sigma))$, that will be specified in the sequel, be arbitrary. For any $n \geq N_0 +1$, $\epsilon \in\{-1,+1\}$, and $t \geq 0$, we introduce
\begin{subequations}
\begin{align}
v_n(t) 
& = c_n(t) - c_n(t-h(t)) , \label{eq: def vn(t)} \\
q_n(t) 
& = \begin{bmatrix} \left< \mathcal{U}L \{ u(t) + d_b(t) \} , \psi_{n,-1} \right> \\ \left< \mathcal{U}L \{ u(t) + d_b(t) \} , \psi_{n,+1} \right> \end{bmatrix} , \label{eq: def qn(t)} \\
r_n(t) 
& = \begin{bmatrix} \left< L \{ u(t) + d_b(t) \} , \psi_{n,-1} \right> \\ \left< L \{ u(t) + d_b(t) \} , \psi_{n,+1} \right> \end{bmatrix} . \label{eq: def rn(t)}
\end{align}
\end{subequations}
Then, based on (\ref{eq: spectral decomposition - prel bis}) and noting that $B_n \{ u(t) + d_b(t) \} = q_n(t) - \Lambda_n r_n(t)$, we have for all $n \geq N_0 +1$ and all $t \geq 0$
\begin{align}
\dot{c}_n(t) 
& = \Lambda_n c_n(t) - M_n v_n(t) \label{eq: modified spectral decomposition for stab analysis inf dim residual dynamics}\\
& \phantom{=}\; + q_n(t) - \Lambda_n r_n(t) + p_{d,n}(t) , \nonumber
\end{align}
with initial condition given over $[-h_M,0]$ by (\ref{eq: spectral decomposition IC - prel bis}). We introduce for any $t \geq -h_M$ the series:
\begin{align*}
S_c(t) & = \sum\limits_{n \geq N_0 + 1} \Vert c_n(t) \Vert^2 = \sum\limits_{\substack{n \geq N_0 + 1 \\ \epsilon\in\{-1,+1\}}} \vert c_{n,\epsilon}(t) \vert^2 
\end{align*}
with, based on (\ref{eq: Riesz basis inequality}), $S_c(t) \leq \Vert \Phi(t) \Vert^2 / m_R$ for all $t \in [-h_M,0]$ while $S_c(t) \leq \Vert X(t) \Vert^2 / m_R$ for all $t \geq 0$. Similarly, we introduce for all $t \geq 0$ the series:
\begin{align*}
S_p(t) & = \sum\limits_{n \geq N_0 + 1} \Vert p_{d,n}(t) \Vert^2 \leq \Vert d_d(t) \Vert^2 / m_R , \\
S_q(t) & = \sum\limits_{n \geq N_0 + 1} \Vert q_n(t) \Vert^2 \leq \dfrac{2 \Vert \mathcal{U}L \Vert^2}{m_R} \left( \Vert u(t) \Vert^2 + \Vert d_b(t) \Vert^2 \right) , \\
S_v(t) & = \sum\limits_{n \geq N_0 + 1} \Vert v_n(t) \Vert^2 \leq 2 \{ S_c(t) + S_c(t-h(t)) \} .
\end{align*}
In particular, using (\ref{eq: assumed estimate on u}) and $\kappa < \sigma$, we obtain that, for all $t \geq 0$,
\begin{align}
S_q(t) 
& \leq \gamma_{1,1} e^{- 2 \kappa t} \Vert \Phi \Vert_{1,h_M}^2 \label{eq: estimate Q} \\
& \phantom{\leq}\; + \gamma_{1,2} \sup\limits_{\tau\in [0,t]} e^{-2\kappa(t-\tau)} \Vert d_d(\tau) \Vert^2 \nonumber \\
& \phantom{\leq}\; + \gamma_{1,3} \sup\limits_{\tau\in [0,t]} e^{-2\kappa(t-\tau)} \Vert d_b(\tau) \Vert^2 \nonumber
\end{align} 
where $\gamma_{1,1} = 6 C_6^2 \Vert \mathcal{U}L \Vert^2/m_R$, $\gamma_{1,2} = 6 C_7^2 \Vert \mathcal{U}L \Vert^2/m_R$, and $\gamma_{1,3} = 2 \Vert \mathcal{U}L \Vert^2 (1+3C_8^2) /m_R$.

We now introduce for any $t_1 < t_2$ and any continuous function $f : [t_1,t_2] \rightarrow \mathbb{R}$ the notation 
\begin{equation}\label{eq: def integral I}
\mathcal{I}(f,t_1,t_2) = \int_{t_1}^{t_2} e^{- 2 \eta (t_2-\tau)} \vert f(\tau) \vert \diff\tau .
\end{equation}
Then, the two following inequalities, which will be useful in the sequel, hold:
\begin{align}
\mathcal{I}(f,t_1,t_2)
& = e^{-2\kappa t_2} \int_{t_1}^{t_2} e^{-2(\eta-\kappa)(t_2-\tau)} \times e^{2\kappa\tau} \vert f(\tau) \vert \diff\tau \nonumber \\
& \leq e^{-2\kappa t_2} \dfrac{1 - e^{-2(\eta-\kappa)(t_2-t_1)}}{2 (\eta-\kappa)} \sup\limits_{\tau\in[t_1,t_2]} e^{2\kappa\tau} \vert f(\tau) \vert \label{eq: integral I - estimate 1}
\end{align}
and, by Cauchy-Schwarz,
\begin{align}
\mathcal{I}(f,t_1,t_2)^2
& \leq \int_{t_1}^{t_2} e^{-2\eta(t_2-\tau)} \diff\tau \times \mathcal{I}(f^2,t_1,t_2) \nonumber \\
& \leq \dfrac{1 - e^{-2\eta(t_2-t_1)}}{2 \eta} \mathcal{I}(f^2,t_1,t_2) . \label{eq: integral I - estimate 2}
\end{align}

For any $t \geq h_M$, we obtain by direct integration of (\ref{eq: modified spectral decomposition for stab analysis inf dim residual dynamics}) over the time interval $[t-h(t),t]$ that
\begin{align*}
v_n(t)
& = (e^{\Lambda_n h(t)}  - I) c_n(t-h(t)) - \Lambda_n \int_{t-h(t)}^{t} e^{\Lambda_n(t-\tau)} r_n(\tau) \diff\tau \\
& \phantom{=}\; + \int_{t-h(t)}^{t} e^{\Lambda_n(t-\tau)} \{ - M_n v_n(\tau) + q_n(\tau) + p_{d,n}(\tau) \} \diff\tau .
\end{align*}
Recalling that $\Lambda_n = \mathrm{diag}(\lambda_{n,-1},\lambda_{n,+1})$ with $\lambda_{n,-1} < \lambda_{n,+1} < 0$ for $n \geq N_0 + 1$ and $h(t) \geq 0$, we have $\Vert e^{\Lambda_n h(t)} - I \Vert \leq 1$ and $\Vert e^{\Lambda_n s} \Vert = e^{\lambda_{n,+1} s} \leq e^{-2\eta s}$ for all $s \geq 0$. Thus, using estimate (\ref{eq: integral I - estimate 2}), as well as estimate (\ref{eq: estimate norm Mn}) derived in Appendix~\ref{annex: estimate norm Mn}, we have for all $t \geq h_M$
\begin{align*}
& \Vert v_n(t) \Vert^2 \\
& \quad \leq 5 \Vert c_n(t-h(t)) \Vert^2 + 5 m_n^2 \mathcal{I}(\Vert v_n \Vert,t-h(t),t)^2 \\
& \quad \phantom{=}\; + 5 \left\Vert \Lambda_n \int_{t-h(t)}^{t} e^{\Lambda_n(t-\tau)} r_n(\tau) \diff\tau \right\Vert^2 \\
& \quad \phantom{=}\; + 5 \mathcal{I}(\Vert q_n \Vert,t-h(t),t)^2 + 5 \mathcal{I}(\Vert p_{d,n} \Vert,t-h(t),t)^2 \\
& \quad \leq 5 \Vert c_n(t-h(t)) \Vert^2 
+ 5 \gamma_2 m_n^2 \mathcal{I}(\Vert v_n \Vert^2,t-h(t),t) \\
& \quad \phantom{=}\; + 5 \left\Vert \Lambda_n \int_{t-h(t)}^{t} e^{\Lambda_n(t-\tau)} r_n(\tau) \diff\tau \right\Vert^2 \\
& \quad \phantom{\leq}\; + 5 \gamma_2 \mathcal{I}(\Vert q_n \Vert^2,t-h(t),t) + 5 \gamma_2 \mathcal{I}(\Vert p_{d,n} \Vert^2,t-h(t),t) 
\end{align*}
with $\gamma_2 = (1-e^{-2\eta h_M})/(2\eta)$. Then, summing for $n \geq N_0 +1$ while using $0 < m_n \leq m_{N_0 + 1}$, we infer that
\begin{align}
S_v(t) & \leq 5 S_c(t-h(t)) + 5 \gamma_2 m_{N_0 + 1}^2 \mathcal{I}(S_v,t-h(t),t) \label{eq: prel estimate Sv} \\
& \phantom{\leq}\; + 5 \sum\limits_{n \geq N_0 + 1} \left\Vert \Lambda_n \int_{t-h(t)}^{t} e^{\Lambda_n(t-\tau)} r_n(\tau) \diff\tau \right\Vert^2 \nonumber \\
& \phantom{\leq}\; + 5 \gamma_2 \mathcal{I}(S_q,t-h(t),t) + 5 \gamma_2 \mathcal{I}(S_p,t-h(t),t) \nonumber
\end{align}
for all $t \geq h_M$. We estimate the third term on the right hand side of the above equation as follows. Let $0 \leq t_0 \leq t$ be given. First, we note that
\begin{align*}
& \Lambda_n \int_{t_0}^{t} e^{\Lambda_n(t-\tau)} r_n(\tau) \diff\tau \\
& \qquad = \begin{bmatrix}
\int_{t_0}^{t} \lambda_{n,-1} e^{\lambda_{n,-1}(t-\tau)} \left< L w(\tau) , \psi_{n,-1} \right> \diff\tau \\
\int_{t_0}^{t} \lambda_{n,+1} e^{\lambda_{n,+1}(t-\tau)} \left< L w(\tau) , \psi_{n,+1} \right> \diff\tau
\end{bmatrix} .
\end{align*}
with $w = u + d_b$. Furthermore, as $\lambda_{n,-1} < \lambda_{n,+1} \leq - 2 \eta < - 2 \kappa < -\kappa < 0$ for $n \geq N_0 + 1$, and introducing $w = (w_1,w_2)$, we have
\begin{align*}
& \left\vert \int_{t_0}^{t} \lambda_{n,\epsilon} e^{\lambda_{n,\epsilon}(t-\tau)} \left< L w(\tau) , \psi_{n,\epsilon} \right> \diff\tau \right\vert \\
& \leq \int_{t_0}^{t} -\lambda_{n,\epsilon} e^{\lambda_{n,\epsilon}(t-\tau)} \vert \left< L w(\tau) , \psi_{n,\epsilon} \right> \vert  \diff\tau \\
& \leq \sum\limits_{m = 1}^2 \int_{t_0}^{t} -\lambda_{n,\epsilon} e^{\lambda_{n,\epsilon}(t-\tau)} \vert \left< L e_m , \psi_{n,\epsilon} \right> \vert \vert w_{m}(\tau) \vert \diff\tau \\
& \leq \sum\limits_{m = 1}^2 \vert \left< L e_m , \psi_{n,\epsilon} \right> \vert \int_{t_0}^{t} -\lambda_{n,\epsilon} e^{\lambda_{n,\epsilon}(t-\tau)} \Vert w(\tau) \Vert \diff\tau \\
& \leq \sum\limits_{m = 1}^2 \vert \left< L e_m , \psi_{n,\epsilon} \right> \vert e^{-\kappa t} \\
& \phantom{\leq}\; \times \int_{t_0}^{t} -\lambda_{n,\epsilon} e^{(\lambda_{n,\epsilon}+\kappa)(t-\tau)} \times e^{\kappa\tau} \Vert w(\tau) \Vert \diff\tau \\
& \leq \dfrac{\lambda_{n,\epsilon}}{\lambda_{n,\epsilon}+\kappa} \sum\limits_{m = 1}^2 \vert \left< L e_m , \psi_{n,\epsilon} \right> \vert e^{-\kappa t} \sup\limits_{\tau\in[t_0,t]} e^{\kappa\tau} \Vert w(\tau) \Vert \\
& \leq \dfrac{2 \eta}{2 \eta - \kappa} \sum\limits_{m = 1}^2 \vert \left< L e_m , \psi_{n,\epsilon} \right> \vert e^{-\kappa t} \sup\limits_{\tau\in[t_0,t]} e^{\kappa\tau} \Vert w(\tau) \Vert .
\end{align*}
Then, we infer
\begin{align*}
& \left\vert \int_{t_0}^{t} \lambda_{n,\epsilon} e^{\lambda_{n,\epsilon}(t-\tau)} \left< L w(\tau) , \psi_{n,\epsilon} \right> \diff\tau \right\vert^2 \\
& \leq \dfrac{8 \eta^2}{(2 \eta - \kappa)^2} \sum\limits_{m = 1}^2 \vert \left< L e_m , \psi_{n,\epsilon} \right> \vert^2 e^{-2\kappa t} \sup\limits_{\tau\in[t_0,t]} e^{2\kappa\tau} \Vert w(\tau) \Vert^2 .
\end{align*}
Then, recalling that $w = u + d_b$ and using (\ref{eq: Riesz basis inequality}) and (\ref{eq: assumed estimate on u}), we obtain:
\begin{align}
& \sum\limits_{n \geq N_0 + 1} \left\Vert \Lambda_n \int_{t_0}^{t} e^{\Lambda_n(t-\tau)} r_n(\tau) \diff\tau \right\Vert^2 \label{eq: estimate series integral rn} \\
& \qquad\leq \gamma_{3,1} e^{-2 \kappa t} \Vert \Phi \Vert_{1,h_M}^2 
+ \gamma_{3,2} e^{-2 \kappa t} \sup\limits_{\tau \in [0,t]} e^{2\kappa\tau} \Vert d_d(\tau) \Vert^2 \nonumber \\
& \qquad\phantom{\leq}\; + \gamma_{3,3} e^{-2 \kappa t} \sup\limits_{\tau \in [0,t]} e^{2\kappa\tau} \Vert d_b(\tau) \Vert^2 \nonumber
\end{align}
with $\gamma_{3,0} = \dfrac{16 \eta^2}{m_R(2 \eta - \kappa)^2} \sum\limits_{m = 1}^2 \Vert L e_m \Vert^2$, $\gamma_{3,1} = 3 \gamma_{3,0} C_6^2$, $\gamma_{3,2} = 3 \gamma_{3,0} C_7^2$, and $\gamma_{3,3} = \gamma_{3,0} (1+3C_8^2)$. Thus, using estimates (\ref{eq: estimate Q}), (\ref{eq: integral I - estimate 1}), and (\ref{eq: estimate series integral rn}) with $t_0 = t-h(t)$ into (\ref{eq: prel estimate Sv}), we infer that
\begin{align*}
S_v(t) & \leq 5 e^{2 \kappa h_M} e^{-2 \kappa h(t)} S_c(t-h(t)) \\
& \phantom{\leq}\; + 5 \gamma_{4,0} m_{N_0 + 1}^2 e^{-2 \kappa t} \sup\limits_{\tau \in [t-h(t),t]} e^{2\kappa\tau} S_v(\tau) \\
& \phantom{\leq}\; + 5 \gamma_{4,1} e^{-2 \kappa t} \Vert \Phi \Vert_{1,h_M}^2 \\
& \phantom{\leq}\; + 5 \gamma_{4,2} e^{-2 \kappa t} \sup\limits_{\tau \in [0,t]} e^{2\kappa\tau} \Vert d_d(\tau) \Vert^2 \\
& \phantom{\leq}\; + 5 \gamma_{4,3} e^{-2 \kappa t} \sup\limits_{\tau \in [0,t]} e^{2\kappa\tau} \Vert d_b(\tau) \Vert^2
\end{align*}
for all $t \geq h_M$ with $\gamma_{4,0} = (1-e^{-2\eta h_M}) (1-e^{-2(\eta-\kappa) h_M})/(4\eta(\eta-\kappa))$, $\gamma_{4,1} = \gamma_{3,1} + \gamma_{1,1}\gamma_{4,0}$, $\gamma_{4,2} = \gamma_{3,2} + \gamma_{4,0} ( \gamma_{1,2} + 1/m_R)$, and $\gamma_{4,3} = \gamma_{3,3} + \gamma_{1,3} \gamma_{4,0}$, thus
\begin{align}
& \sup\limits_{\tau\in[h_M,t]} e^{2\kappa\tau} S_v(\tau) \nonumber \\
& \quad \leq 5 e^{2\kappa h_M} \sup\limits_{\tau\in[0,t]} e^{2\kappa\tau} S_c(\tau) \label{eq: prel-prel estimate sup V} \\
& \quad \phantom{\leq}\; + 5 \gamma_{4,0} m_{N_0 + 1}^2 \sup\limits_{\tau\in[0,t]} e^{2\kappa\tau} S_v(\tau) 
+ 5 \gamma_{4,1} \Vert \Phi \Vert_{1,h_M}^2 \nonumber \\
& \quad \phantom{\leq}\; + 5 \gamma_{4,2} \sup\limits_{\tau\in[0,t]} e^{2 \kappa \tau} \Vert d_d(\tau) \Vert^2 
+ 5 \gamma_{4,3} \sup\limits_{\tau\in[0,t]} e^{2 \kappa \tau} \Vert d_b(\tau) \Vert^2 . \nonumber
\end{align}

In order to estimate the $S_c$-term in (\ref{eq: prel-prel estimate sup V}), we integrate (\ref{eq: modified spectral decomposition for stab analysis inf dim residual dynamics}) over the time interval $[0,t]$ for $t \geq 0$, yielding the following estimate:
\begin{align*}
\Vert c_n(t) \Vert & \leq e^{-2\eta t} \Vert c_n(0) \Vert + m_{N_0+1} \mathcal{I}(\Vert v_n \Vert,0,t) \\
& \phantom{\leq}\; + \left\Vert \Lambda_n \int_{0}^{t} e^{\Lambda_n(t-\tau)} r_n(\tau) \diff\tau \right\Vert \\ 
& \phantom{\leq}\; + \mathcal{I}(\Vert q_n \Vert,0,t) + \mathcal{I}(\Vert p_{d,n} \Vert,0,t) ,
\end{align*}
where we used estimate (\ref{eq: estimate norm Mn}) derived in Appendix~\ref{annex: estimate norm Mn} and the fact that $\lambda_{n,\epsilon} \leq - 2 \eta$. Using (\ref{eq: integral I - estimate 2}), we deduce that
\begin{align*}
\Vert c_n(t) \Vert^2 & \leq 5 e^{-4\eta t} \Vert c_n(0) \Vert^2 + \dfrac{5}{2\eta} m_{N_0+1}^2 \mathcal{I}(\Vert v_n \Vert^2,0,t) \\
& \phantom{\leq}\; + 5 \left\Vert \Lambda_n \int_{0}^{t} e^{\Lambda_n(t-\tau)} r_n(\tau) \diff\tau \right\Vert^2 \\ 
& \phantom{\leq}\; + \dfrac{5}{2\eta} \mathcal{I}(\Vert q_n \Vert^2,0,t) + \dfrac{5}{2\eta} \mathcal{I}(\Vert p_{d,n} \Vert^2,0,t)
\end{align*}
and thus, using (\ref{eq: integral I - estimate 1}), (\ref{eq: estimate series integral rn}) with $t_0 = 0$, and $S_c(0) \leq \Vert \Phi(0) \Vert^2 / m_R$,
\begin{align*}
S_c(t)
& \leq 5 e^{-2\kappa t} S_c(0) + \dfrac{5}{2\eta} m_{N_0+1}^2 \mathcal{I}(S_v,0,t) \\
& \phantom{\leq}\; + 5 \sum\limits_{n \geq N_0+1} \left\Vert \Lambda_n \int_{0}^{t} e^{\Lambda_n(t-\tau)} r_n(\tau) \diff\tau \right\Vert^2 \\
& \phantom{\leq}\; + \dfrac{5}{2\eta} \mathcal{I}(S_q,0,t) + \dfrac{5}{2\eta} \mathcal{I}(S_p,0,t) \\
& \leq 5 \gamma_{5,0} m_{N_0+1}^2 e^{-2 \kappa t} \sup\limits_{\tau \in [0,t]} e^{2\kappa\tau} S_v(\tau) \\
& \phantom{\leq}\; + 5 \gamma_{5,1} e^{-2 \kappa t} \Vert \Phi \Vert_{1,h_M}^2 \\
& \phantom{\leq}\; + 5 \gamma_{5,2} e^{-2\kappa t} \sup\limits_{\tau \in [0,t]} e^{2\kappa\tau} \Vert d_d(\tau) \Vert^2 \\
& \phantom{\leq}\; + 5 \gamma_{5,3} e^{-2\kappa t} \sup\limits_{\tau \in [0,t]} e^{2\kappa\tau} \Vert d_b(\tau) \Vert^2
\end{align*}
for all $t \geq 0$, with $\gamma_{5,0} = 1/(4\eta(\eta-\kappa))$, $\gamma_{5,1} = \gamma_{3,1} + \gamma_{1,1}\gamma_{5,0} + 1/m_R$, $\gamma_{5,2} = \gamma_{3,2} + \gamma_{5,0} ( \gamma_{1,2} + 1/m_R)$, and $\gamma_{5,3} = \gamma_{3,3} + \gamma_{1,3} \gamma_{5,0}$. This yields the estimate
\begin{align}
& \sup\limits_{\tau \in [0,t]} e^{2\kappa\tau} S_c(\tau) \nonumber \\
& \quad \leq 5 \gamma_{5,0} m_{N_0+1}^2 \sup\limits_{\tau \in [0,t]} e^{2\kappa\tau} S_v(\tau) \label{eq: estimate sup Z over [hm,t]} \\
& \quad \phantom{\leq}\; + 5 \gamma_{5,1} \Vert \Phi \Vert_{1,h_M}^2 \nonumber 
+ 5 \gamma_{5,2} \sup\limits_{\tau \in [0,t]} e^{2\kappa\tau} \Vert d_d(\tau) \Vert^2 \nonumber \\
& \quad \phantom{\leq}\; + 5 \gamma_{5,3} \sup\limits_{\tau \in [0,t]} e^{2\kappa\tau} \Vert d_b(\tau) \Vert^2 \nonumber
\end{align}
for all $t \geq 0$. Now, introducing the quantity:
\begin{align*}
\delta & =  5 m_{N_0+1}^2 ( \gamma_{4,0} + 5 e^{2 \kappa h_M} \gamma_{5,0} ) \\
& = \dfrac{10 \alpha^2 \gamma^2}{(\alpha^2-1) (N_0+1)^4 \pi^4 + \beta} \\
& \phantom{=}\; \times \dfrac{(1-e^{-2\eta h_M}) (1-e^{-2(\eta-\kappa) h_M}) + 5 e^{2 \kappa h_M}}{4\eta(\eta-\kappa)} ,
\end{align*}
the combination of (\ref{eq: prel-prel estimate sup V}-\ref{eq: estimate sup Z over [hm,t]}) yields, for all $t \geq h_M$, 
\begin{align*}
& \sup\limits_{\tau\in[h_M,t]} e^{2\kappa\tau} S_v(\tau) \\
& \leq \delta \sup\limits_{\tau\in[0,t]} e^{2\kappa\tau} S_v(\tau) + \gamma_{6,1} \Vert \Phi \Vert_{1,h_M}^2 \nonumber \\
& \phantom{\leq}\; + \gamma_{6,2} \sup\limits_{\tau\in[0,t]} e^{2 \kappa \tau} \Vert d_d(\tau) \Vert^2 
+ \gamma_{6,3} \sup\limits_{\tau\in[0,t]} e^{2 \kappa \tau} \Vert d_b(\tau) \Vert^2 , \nonumber
\end{align*}
where $\gamma_{6,i} = 5 (\gamma_{4,i} + 5 e^{2\kappa h_M} \gamma_{5,i})$. From the small gain condition (\ref{eq: thm infinite dim residual - small gain condition}), a continuity argument in $\kappa = 0$ shows the existence of a constant $\kappa \in (0,\min(\eta,\sigma))$, independent of $\Phi, d_d, d_b, h, u$, such that $\delta < 1$. Selecting such a $\kappa$ for the remaining of the proof, we obtain that
\begin{align}
& \sup\limits_{\tau\in[h_M,t]} e^{2\kappa\tau} S_v(\tau) \nonumber \\
& \leq \dfrac{\delta}{1-\delta} \sup\limits_{\tau\in[0,h_M]} e^{2\kappa\tau} S_v(\tau) + \dfrac{\gamma_{6,1}}{1-\delta} \Vert \Phi \Vert_{1,h_M}^2 \label{eq: estimate sup V over [2hM,t]} \\
& \phantom{\leq}\; + \dfrac{\gamma_{6,2}}{1-\delta} \sup\limits_{\tau\in[0,t]} e^{2 \kappa \tau} \Vert d_d(\tau) \Vert^2 
+ \dfrac{\gamma_{6,3}}{1-\delta} \sup\limits_{\tau\in[0,t]} e^{2 \kappa \tau} \Vert d_b(\tau) \Vert^2 \nonumber
\end{align}
for all $t \geq h_M$. 

Now, based on (\ref{eq: modified spectral decomposition for stab analysis inf dim residual dynamics}), direct estimations reported in Appendix~\ref{annex: rough estimate Sc over [0,hm]} show the existence of constants $\gamma_{7,0},\gamma_{7,1},\gamma_{7,2},\gamma_{7,3} > 0$, independent of $\Phi, d_d, d_b, h, u$, such that $S_c(\tau) \leq \gamma_{7,0} \Vert \Phi \Vert_{1,h_M}^2$ for all $\tau \in [-h_M,0]$ and 
\begin{align}
S_c(t) & \leq \gamma_{7,1} \Vert \Phi \Vert_{1,h_M}^2 + \gamma_{7,2} \sup\limits_{\tau \in [0,t]} e^{2\kappa\tau} \Vert d_d(\tau) \Vert^2 \label{eq: rough estimate Sc over [0,hm]}\\
& \phantom{\leq}\; + \gamma_{7,3} \sup\limits_{\tau \in [0,t]} e^{2\kappa\tau} \Vert d_b(\tau) \Vert^2 \nonumber
\end{align}
for all $0 \leq t \leq h_M$. Consequently, we have for all $0 \leq t \leq h_M$
\begin{align}
\sup\limits_{\tau\in[0,t]} e^{2\kappa\tau} S_v(\tau)
& \leq 4 e^{2 \kappa h_M} \sup\limits_{\tau\in[-h_M,t]} S_c(\tau) \nonumber \\
& \leq 4 (\gamma_{7,0} + \gamma_{7,1}) e^{2 \kappa h_M} \Vert \Phi \Vert_{1,h_M}^2 \label{eq: estimate sup V over [hM,t] for t leq 2hm} \\
& \phantom{\leq}\; + 4 \gamma_{7,2} e^{2 \kappa h_M} \sup\limits_{\tau \in [0,t]} e^{2\kappa\tau} \Vert d_d(\tau) \Vert^2 \nonumber \\
& \phantom{\leq}\; + 4 \gamma_{7,3} e^{2 \kappa h_M} \sup\limits_{\tau \in [0,t]} e^{2\kappa\tau} \Vert d_b(\tau) \Vert^2 \nonumber .
\end{align}
From estimates (\ref{eq: estimate sup V over [2hM,t]}) and (\ref{eq: estimate sup V over [hM,t] for t leq 2hm}), we infer the existence of constants $\gamma_{8,1},\gamma_{8,2},\gamma_{8,3} > 0$, independent of $\Phi, d_d, d_b, h, u$, such that
\begin{align}
\sup\limits_{\tau\in[0,t]} e^{2\kappa\tau} S_v(\tau)
& \leq \gamma_{8,1} \Vert \Phi \Vert_{1,h_M}^2 + \gamma_{8,2} \sup\limits_{\tau \in [0,t]} e^{2\kappa\tau} \Vert d_d(\tau) \Vert^2 \label{eq: estimate sup V over [hM,t]} \\
& \phantom{\leq}\; + \gamma_{8,3} \sup\limits_{\tau \in [0,t]} e^{2\kappa\tau} \Vert d_b(\tau) \Vert^2 \nonumber
\end{align}
for all $t \geq 0$. Substituting (\ref{eq: estimate sup V over [hM,t]}) into (\ref{eq: estimate sup Z over [hm,t]}), we obtain the existence of constants $C_9,C_{10}, C_{11} > 0$, independent of $\Phi,d_d,d_b,h,u$, such that the claimed estimate (\ref{eq: exp stab infinite-dim part}) holds.
\qed

\subsection{Proof of the main result}

We are now ready to prove the main result of this paper.

\textbf{Proof of Theorem~\ref{thm: main theorem}.}
We select $N_0 \geq 1$ large enough such that $N_0 \geq \lfloor \beta^{1/4}/\pi \rfloor$ and (\ref{eq: main thm - small gain condition}) holds. With the matrices $A$ and $B$ defined by (\ref{eq: truncated model - A}-\ref{eq: truncated model - B}), Lemma~\ref{lem: kalman condition} ensures the existence of a feedback gain $K \in \mathbb{R}^{2 \times 2 N_0}$ such that $A_\mathrm{cl} = A + BK$ is Hurwitz. Then, based on Lemma~\ref{lem: prel lemma feasability}, let $0 < h_m < h_M$ be such that the LMI (\ref{eq: main thm - LMI condition}) is feasible. By a continuity argument, there exists $\sigma > 0$ such that $\Theta(h_M,\sigma) \prec 0$ where $\Theta(h_M,\sigma)$ is defined by (\ref{eq: prel lemma - LMI condition}) with $F = A_\mathrm{cl}$ and $G = M$ given by (\ref{eq: truncated model - M}). Then, Lemma~\ref{lem: stab truncated model} provides constants $C_3, C_4, C_5 > 0$ such that, for any distributed perturbation $d_d \in \mathcal{C}^0(\mathbb{R}_+;L^2(0,1))$, any boundary perturbation $d_b \in \mathcal{C}^1(\mathbb{R}_+;\mathbb{R}^2)$, any initial condition $\Phi \in \mathcal{C}^1([-h_M,0];\mathcal{H})$, and any delay $h \in \mathcal{C}^0(\mathbb{R}_+;\mathbb{R})$ with $h_m \leq h \leq h_M$, the mild solution $X\in\mathcal{C}^0(\mathbb{R}_+;\mathcal{H})$ of (\ref{eq: abstract form - DE}-\ref{eq: abstract form - IC}) with $u = K Y \in \mathcal{C}^1(\mathbb{R}_+;\mathbb{R}^2)$ satisfies 
\begin{align}
\sqrt{ \sum\limits_{\substack{1 \leq n \leq N_0 \\ \epsilon\in\{-1,+1\}}} \vert c_{n,\epsilon}(t) \vert^2 }  
& \leq \dfrac{C_3}{\sqrt{m_R}} e^{- \sigma t} \Vert \Phi \Vert_{1,h_M} \label{eq: proof main result - 1} \\
& \phantom{\leq}\; + \dfrac{C_4}{\sqrt{m_R}} \sup\limits_{\tau \in [0,t]} e^{- \sigma(t-\tau)} \Vert d_d(\tau) \Vert \nonumber \\
& \phantom{\leq}\; + C_5 \sup\limits_{\tau \in [0,t]} e^{- \sigma(t-\tau)} \Vert d_b(\tau) \Vert \nonumber
\end{align}
for all $t \geq 0$, where we recall that $m_R,M_R$ are constants associated with the Riesz basis $\mathcal{F}_\phi$, see Definition~\ref{def: Riesz basis} and Lemma~\ref{lem: Riesz basis}. Now, with $u = KY$, we infer the existence of constants $C_6 , C_7 , C_8 > 0$, independent of $\Phi, d_d, d_b, h$, such that the estimate (\ref{eq: assumed estimate on u}) holds. Applying Lemma~\ref{lem: infinite-dim part negelected in the design}, we have constants $\kappa\in(0,\sigma)$ and $C_9 , C_{10} , C_{11} > 0$, independent of $\Phi, d_d, d_b, h$, such that the estimate (\ref{eq: exp stab infinite-dim part}) holds. Putting together (\ref{eq: exp stab infinite-dim part}) and (\ref{eq: proof main result - 1}), and using the Riesz basis inequality (\ref{eq: Riesz basis inequality}), 
\begin{align*}
\Vert X(t) \Vert
& \leq \sqrt{ M_R \sum\limits_{\substack{n \leq 1 \\ \epsilon\in\{-1,+1\}}} \vert c_{n,\epsilon}(t) \vert^2 } \\
& \leq C_0 e^{- \kappa t} \Vert \Phi(\tau) \Vert_{1,h_M} + C_1 \sup\limits_{\tau \in [0,t]} e^{- \kappa(t-\tau)} \Vert d_d(\tau) \Vert , \nonumber \\
& \phantom{\leq}\; + C_2 \sup\limits_{\tau \in [0,t]} e^{- \kappa(t-\tau)} \Vert d_b(\tau) \Vert , \nonumber
\end{align*}
with $C_0 = \sqrt{M_R} ( C_3/\sqrt{m_R} + \sqrt{C_9} )$, $C_1 = \sqrt{M_R} ( C_4/\sqrt{m_R} + \sqrt{C_{10}} )$, and $C_2 = \sqrt{M_R} ( C_5 + \sqrt{C_{11}} )$ which are constants independent of $\Phi, d_d, d_b, h$. This completes the proof.  
\qed

\section{Numerical example}\label{sec: numerical illustration}
In this section, we describe a numerical example that illustrates the main result of this paper, namely: Theorem~\ref{thm: main theorem}. For numerical computations, we set $\alpha = 1.5$, $\beta_0 = 50$, and $\gamma = 50$. For this choice of parameters, we have one unstable eigenvalue for $\mathcal{U}_0$. Selecting $N_0 = 2$, the small gain condition (\ref{eq: main thm - small gain condition}) is satisfied. Depending on the selected actuation configuration, i.e. either one or two control inputs (see Remark~\ref{rem: number of control input}), we compute first a feedback gain $K \in \mathbb{R}^{2 \times 2N_0}$ such that the eigenvalues of $A+BK$ are given by $\{-5,-6,-7,-8\}$. Then, we compute the corresponding value $h_M > 0$ of the maximum amplitude of the admissible state-delay $h$ provided by the LMI (\ref{eq: main thm - LMI condition}). Obviously, the obtained value of $h_M$ depends on the selected actuation scheme. The numerical computations yield the results reported hereafter\footnote{The provided values $h$ of $h_M$ are such that the LMI (\ref{eq: main thm - LMI condition}) is feasible for $h_M = h$ while no feasible solution was found for $h_M = h + 0.001$ using \textsc{Matlab}~2017b LMI solvers.}. Note that, due to the invariance of the beam model (\ref{eq: beam model - PDE}-\ref{eq: beam model - IC2}) under the change of variable $x \leftarrow 1-x$, the numerical results obtained for the two single input configurations are equivalent.
\begin{itemize}
\item One control input located at $x = 0$:
\begin{equation*}
K = \begin{bmatrix}
   -1.7614    & 0.0276   & 11.8714   & -0.0360 \\
         0    &      0   &       0   &       0
\end{bmatrix}
\end{equation*}
with $h_M = 0.038$.
\item Two control inputs located respectively at $x = 0$ and $x = 1$. 
\begin{equation*}
K = \begin{bmatrix}
    2.0076    & 0.4186   &  5.0313   &  0.1129 \\
    1.9972    & 0.4278   & -4.5575   & -0.0178
\end{bmatrix}
\end{equation*}
with $h_M = 0.239$.
\end{itemize}
We observe that the use of two control inputs allows a significantly larger value of $h_M$. This indicates that, in our numerical setting, the two control inputs configuration is much more robust with respect to state-delays than the single input configuration.

Note that the above values of $h_M$ yield the satisfaction of the ISS-like estimate (\ref{eq: main theorem - estimate y}). However, as discussed in Remark~\ref{rem: ISS-like estimate}, one can get the ISS estimate (\ref{eq: main theorem - ISS in strict form - estimate y}) when considering the small gain condition (\ref{eq: main thm - ISS in strict form}) instead of the LMI condition (\ref{eq: main thm - LMI condition}). However, this yields much more conservative results since, based on the upper estimate (\ref{eq: main thm - ISS in strict form - upperbound}), we must have $h_M < 0.0070$ in the single input configuration while $h_M < 0.0148$ in the two inputs configuration.

For numerical illustration, we consider the two control inputs setting and we set the parameters of the simulation as follows: initial conditions $y_0(t,x) = 2 (1-t)^2 x (1-x)$ and $y_{t0}(t,x) = - (1-t)^2 \sin(4\pi x)(1+2x)$; time varying delay $h(t) = 0.12 + 0.1 \sin(6 \pi t)$; distributed perturbation $d_d(t,x) = 3 \exp(- 2 (t-5)^2) ( 2 + \cos(2 \pi x))$; and boundary perturbations $d_{b,1}(t) = \cos(2\pi t)\exp(-2(t-5)^2)$ and $d_{b,2}(t) = -\sin(3\pi t)\exp(-2(t-5)^2)$. In particular, the perturbations are vanishing ones with maximum amplitudes at time $t=5\,\mathrm{s}$. The used numerical scheme consists in the modal approximation of the Euler-Bernouilli beam model (\ref{eq: beam model - PDE}-\ref{eq: beam model - IC2}) using its first 40 modes. 

\begin{figure}
    \centering
	\includegraphics[width=3.5in]{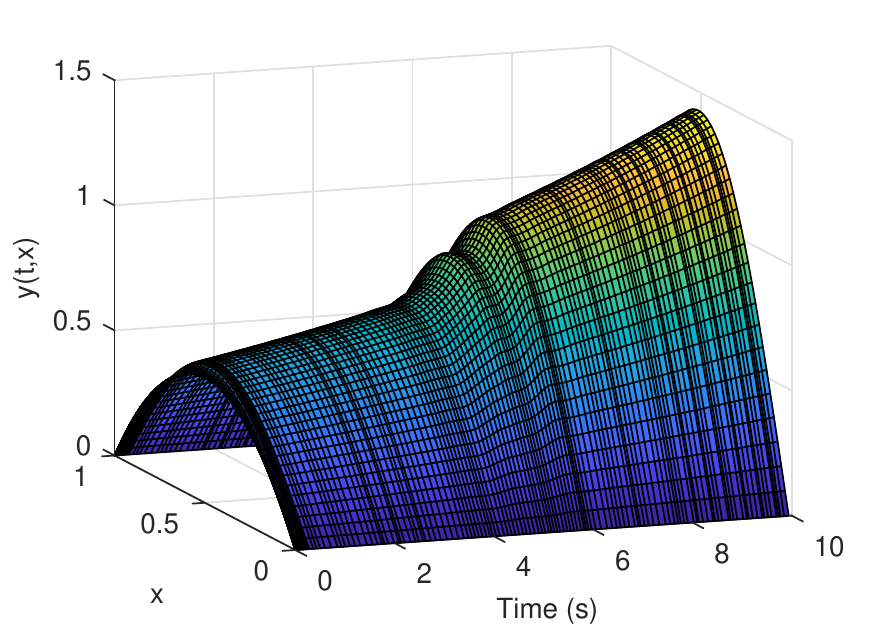}
    \caption{Time evolution of the displacement $y(t,x)$ for the open-loop system}
    \label{fig: sim - open-loop}
\end{figure}

\begin{figure}
     \centering
     	\subfigure[Displacement $y(t,x)$]{
		\includegraphics[width=3.5in]{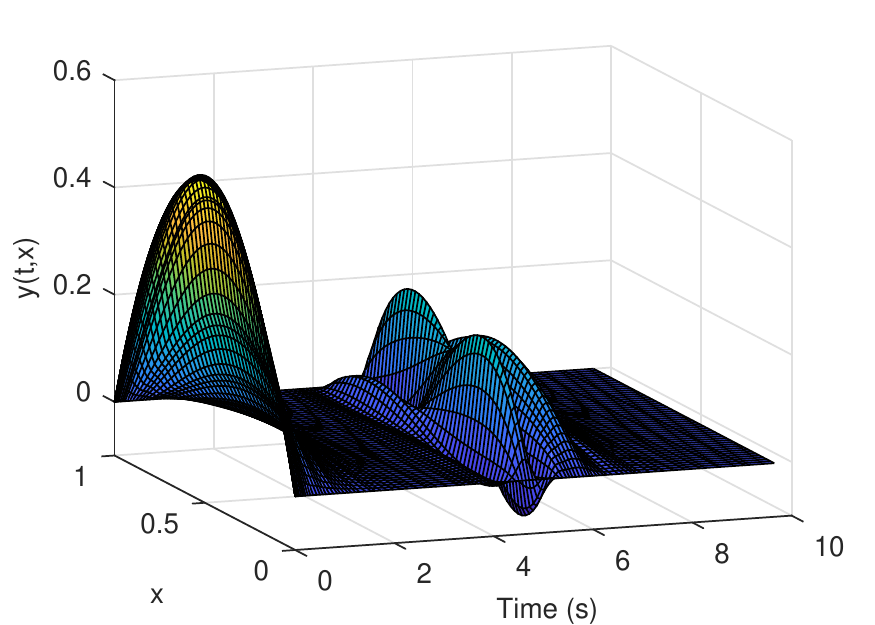}
		\label{fig: sim closed-loop - displacement}
		}
		\subfigure[Speed of displacement $y_t(t,x)$]{
		\includegraphics[width=3.5in]{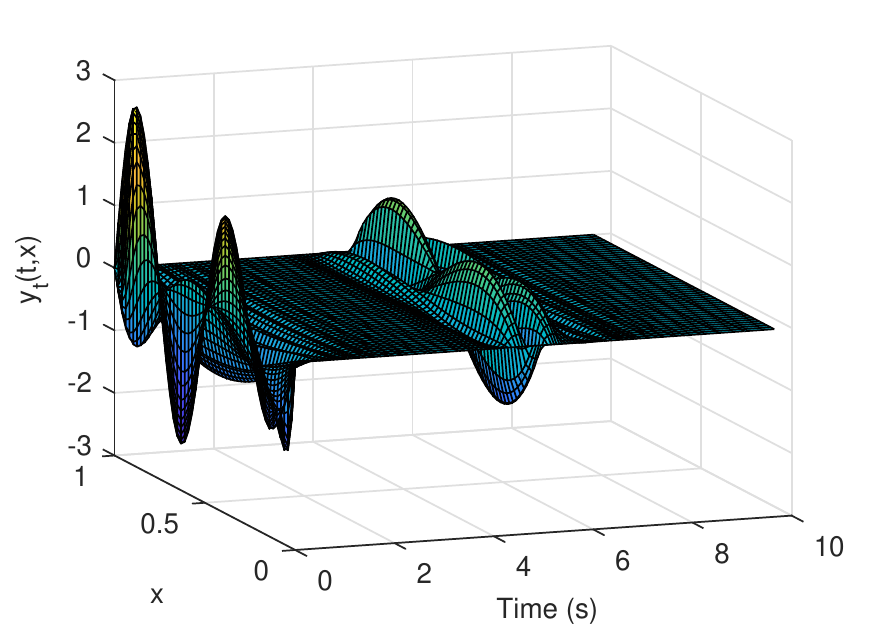}
		\label{fig: sim closed-loop - speed}
		}
     	\subfigure[Command inputs $u_1$ and $u_2$]{
		\includegraphics[width=3.5in]{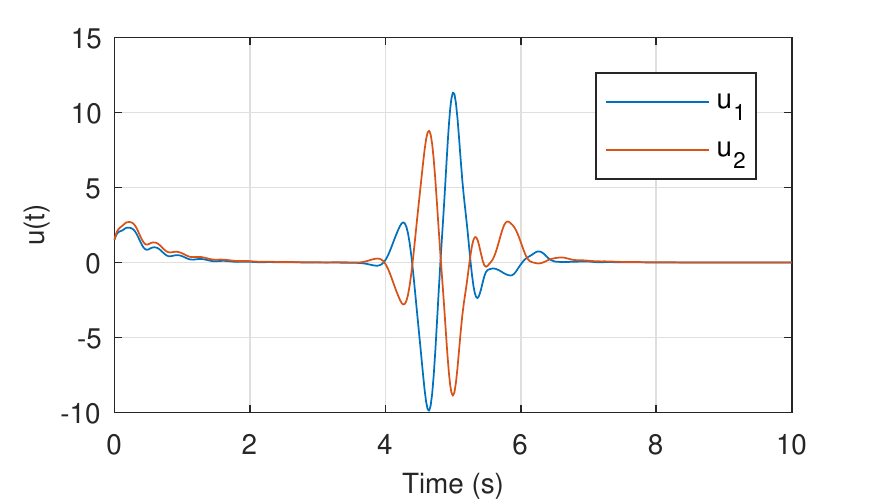}
		\label{fig: sim closed-loop - control}
		}
     \caption{Time evolution of the closed-loop system with two command inputs $u_1$ and $u_2$}
     \label{fig: sim - closed-loop}
\end{figure}

The time domain evolution of the open-loop system is depicted in Fig.~\ref{fig: sim - open-loop}, illustrating the unstable behavior of the studied Euler-Bernoulli beam. In this configuration, numerical computations show that both single input configurations fail to ensure the stability of the Euler-Bernoulli beam in closed-loop. This is due to the fact that the maximum amplitude of the considered input delay is $0.22$, which is significantly higher that the guaranteed bound $h_M=0.037$ obtained via the LMI condition. Conversely, the temporal behavior of the closed-loop system in the two control inputs configuration is depicted in Fig~\ref{fig: sim - closed-loop}. We observe that the proposed control strategy achieves the exponential stabilization of the studied Euler-Bernoulli beam by quickly damping out the initial condition with control inputs converging to zero. Furthermore, the impact of the vanishing distributed and boundary perturbations, whose maximum intensity occur at time $t = 5\,\mathrm{s}$, is rapidly eliminated as $t$ increases. These numerical results are in conformity with the fading memory nature of the estimates (\ref{eq: main theorem - estimate y}) and (\ref{eq: cor main theorem - estimate displacement sup norm}).

\section{Conclusion}\label{sec: conclusion}
This paper discussed the point torque boundary feedback stabilization of a damped Euler-Bernoulli beam in the presence of a state-delay. The proposed control law takes the form of a state feedback computed based on a finite-dimensional truncated model of the original infinite-dimensional system. In order to ensure the exponential stability of the closed-loop truncated model in the presence of the state-delay, an LMI-based sufficient condition was derived on the maximum amplitude of the state-delay. Then, provided the fact that the truncated model captures a sufficiently large number of modes (including the unstable ones), the stability of the closed-loop infinite-dimensional system was assessed by resorting to a small gain argument. In the presence of both distributed and boundary disturbances, the derived stability result takes the form of an ISS-like estimate with fading memory of the perturbations. The validity of the proposed control strategy was illustrated with numerical simulations.



\appendices

\section{Proof of (\ref{eq: def mild solution bis})}\label{annex: mild solution}

Introducing $D = \gamma\pi \in \mathcal{L}(\mathcal{H})$, it is a general result that~\cite[Thm.~3.2.1]{Curtain2012}:
\begin{equation}\label{eq: connexion between S and T}
T(t)z = S(t)z + \int_0^t T(t-s) D S(s) z \diff s  \end{equation}
for all $z\in\mathcal{H}$ and all $t \geq 0$. Introducing $a_0 = \Phi(0) - Lw(0)$, $f(t) = \mathcal{A}Lw(t)-L\dot{w}(t)+DX(t-h(t))+p_d(t)$, and $g(t) = f(t) - DX(t)$, the use of (\ref{eq: connexion between S and T}) into (\ref{eq: def mild solution}) yields
\begin{align*}
X(t) 
& = S(t) a_0 + L w(t) + \int_0^t S(t-s) f(s) \diff s \\
& = T(t) a_0 + L w(t) + \int_0^t T(t-s) g(s) \diff s \\
& \phantom{=}\, + \int_0^t T(t-s)D \{ - S(s)a_0 + X(s) \} \diff s - \mathcal{J}(t)
\end{align*}
with
\begin{align*}
\mathcal{J}(t)
& = \int_0^t \int_0^{t-s} T(t-s-\tau) D S(\tau) f(s) \diff\tau\diff s \\
& = \int_0^t \int_0^{t-s} T(\tau) D S(t-s-\tau) f(s) \diff\tau\diff s \\
& = \int_0^t \int_0^{t-\tau} T(\tau) D S(t-s-\tau) f(s) \diff s \diff\tau \\
& = \int_0^t T(\tau) D \int_0^{t-\tau} S(t-s-\tau) f(s) \diff s \diff\tau \\
& = \int_0^t T(t-\xi) D \int_0^{\xi} S(\xi-s) f(s) \diff s \diff\xi
\end{align*}
where two changes of variables and Fubini's theorem have been used. Combining the above identities, we obtain that
\begin{align*}
X(t) 
& = T(t) a_0 + L w(t) + \int_0^t T(t-s) \{ g(s) + DLw(s) \} \diff s .
\end{align*}
Recalling that $\mathcal{U}=\mathcal{A}+D$ with $D = \gamma\Pi$, we obtain the identity (\ref{eq: def mild solution bis}). Similarly, using $T(t)z = S(t)z + \int_0^t S(t-s) D T(s) z \diff s$ (see~\cite[Thm.~3.2.1]{Curtain2012}), we have that (\ref{eq: def mild solution bis}) implies  (\ref{eq: def mild solution}).

\section{Proof of Lemma~\ref{lem: Riesz basis}}\label{annex: Riesz basis}

We first show that $\mathcal{F}_\phi$ is maximal in $\mathcal{H}$. To do so, let $z = (z_1,z_2) \in \mathcal{H}$ be such that $\left< z, \phi_{n,\epsilon} \right> = 0$ for all $n \geq 1$ and $\epsilon \in \{-1,+1\}$. Then we have
\begin{equation*}
\begin{bmatrix}
- n^2 \pi^2 & \lambda_{n,-1} \\ - n^2 \pi^2 & \lambda_{n,+1}
\end{bmatrix}
\begin{bmatrix}
\left[ z''_1 , \sin(n\pi\cdot) \right] \\ 
\left[ z_2 , \sin(n\pi\cdot) \right]
\end{bmatrix}
= 0 ,
.
\end{equation*}
As $\lambda_{n,-1} \neq \lambda_{n,+1}$, we deduce that $\left[ z''_1 , \sin(n\pi\cdot) \right] = \left[ z_2 , \sin(n\pi\cdot) \right] = 0$ for all $n \geq 1$. As $\{ \sin(n\pi\cdot) \,:\, n \geq 1 \}$ is maximal in $L^2(0,1)$, we infer that $z''_1=z_2=0$. From $z_1 \in H_0^1(0,1)$, we deduce that $z=0$ showing that $\mathcal{F}_\phi$ is maximal in $\mathcal{H}$.

We now show that $\mathcal{F}_\phi$ satisfies a Riesz basis inequality of the type (\ref{eq: def Riesz basis inequality}). First, it is straightforward to note that $\left< \phi_{n_1,\epsilon_1} , \phi_{n_2,\epsilon_2} \right> = 0$ for any $n_1 \neq n_2$ and any $\epsilon_1,\epsilon_2\in\{-1,+1\}$. Let $N \geq 1$ and $a_{n,\epsilon} \in \mathbb{C}$ be arbitrary. We have
\begin{align*}
& \left\Vert \sum\limits_{n = 1}^N a_{n,-1} \phi_{n,-1} + a_{n,+1} \phi_{n,+1} \right\Vert^2 \\
& = \sum\limits_{n = 1}^N \left\Vert a_{n,-1} \phi_{n,-1} + a_{n,+1} \phi_{n,+1} \right\Vert^2 \\
& = \sum\limits_{n = 1}^N \vert a_{n,-1} \vert^2 + \vert a_{n,+1} \vert^2 + 2 \operatorname{Re}\left( a_{n,-1} \overline{a_{n,+1}} \left< \phi_{n,-1} , \phi_{n,+1} \right> \right)
\end{align*}
with
\begin{equation*}
\left< \phi_{n,-1} , \phi_{n,+1} \right>
= \dfrac{2 n^4 \pi^4 - \beta}{2 k_{n,-1} k_{n,+1}}
= \dfrac{n^4 \pi^4 - \beta/2}{\sqrt{\alpha^2 n^8 \pi^8 + \beta^2/4}} .
\end{equation*}
As $\varphi : \mathbb{R}_+ \rightarrow (x-\beta/2)/\sqrt{\alpha^2 x^2 + \beta^2/4}$ is strictly increasing with $\varphi(x) \rightarrow 1/\alpha < 1$ when $x \rightarrow + \infty$ and $\varphi(\pi^4) \geq -\beta/\sqrt{4 \alpha^2 \pi^8 + \beta^2} > - 1$, we have $\vert \left< \phi_{n,-1} , \phi_{n,+1} \right> \vert \leq C_R = \max(1/\alpha , \beta/\sqrt{4 \alpha^2 \pi^8 + \beta^2}) < 1$ for all $n \geq 1$. This yields
\begin{align*}
& \vert a_{n,-1} \vert^2 + \vert a_{n,+1} \vert^2 + 2 \operatorname{Re}\left( a_{n,-1} \overline{a_{n,+1}} \left< \phi_{n,-1} , \phi_{n,+1} \right> \right) \\
& \quad\leq \vert a_{n,-1} \vert^2 + \vert a_{n,+1} \vert^2 + 2 C_R \vert a_{n,-1} \vert \vert a_{n,+1} \vert \\
& \quad\leq (1+C_R) \left( \vert a_{n,-1} \vert^2 + \vert a_{n,+1} \vert^2 \right)
\end{align*}
and 
\begin{align*}
& \vert a_{n,-1} \vert^2 + \vert a_{n,+1} \vert^2 + 2 \operatorname{Re}\left( a_{n,-1} \overline{a_{n,+1}} \left< \phi_{n,-1} , \phi_{n,+1} \right> \right) \\
& \quad\geq \vert a_{n,-1} \vert^2 + \vert a_{n,+1} \vert^2 - 2 C_R \vert a_{n,-1} \vert \vert a_{n,+1} \vert \\
& \quad\geq (1-C_R) \left( \vert a_{n,-1} \vert^2 + \vert a_{n,+1} \vert^2 \right).
\end{align*}
The claimed Riesz basis inequality holds with $m_R = 1-C_R \in(0,1)$ and $M_R = 1+C_R > 1$.

\section{Well-posedness of the abstract system (\ref{eq: abstract form - DE}-\ref{eq: abstract form - IC}) with $u = KY$}\label{annex: well-posedness}

Let $d_d \in \mathcal{C}^0(\mathbb{R}_+;L^2(0,1))$, $d_b \in \mathcal{C}^1(\mathbb{R}_+;\mathbb{R}^2)$, $\Phi \in \mathcal{C}^0([-h_M,0];\mathcal{H})$, and $h \in \mathcal{C}^0(\mathbb{R}_+;\mathbb{R})$ with $0<h_m \leq h \leq h_M$ be given. First, we show that $X\in\mathcal{C}^0(\mathbb{R}_+;\mathcal{H})$ such that $Y \in \mathcal{C}^1(\mathbb{R}_+;\mathbb{R}^{2 N_0})$ is a mild solution of (\ref{eq: abstract form - DE}-\ref{eq: abstract form - IC}) with $u = K Y$ if and only if $X\in\mathcal{C}^0(\mathbb{R}_+;\mathcal{H})$ and satisfies for all $t \geq 0$
\begin{align}
& X(t) 
= S(t) \{ \Phi(0) - L \tilde{w}(0) \} + L \tilde{w}(t) \label{eq: def aux mild solution} \\
& \, + \int_0^t S(t-s) \{ \mathcal{A}L\tilde{w}(s) - L\dot{\tilde{w}}(s) + \gamma \Pi X(s-h(s)) + p_d(s) \} \diff s \nonumber
\end{align}
with $\tilde{w} = Kv + d_b$,
\begin{align}
v(t) & = e^{(A_\mathrm{cl}-M) t} Y_\phi(0) \label{eq: aux abstract form - aux input v} \\
& \phantom{=}\; + \int_0^t e^{(A_\mathrm{cl}-M)(t-\tau)} M Y(\tau-h(\tau)) \diff\tau \nonumber \\
& \phantom{=}\; + \int_0^t e^{(A_\mathrm{cl}-M)(t-\tau)} \{B d_b(\tau) + P_d(\tau) \} \diff\tau , \nonumber
\end{align}
the initial condition $X(\tau) = \Phi(\tau)$ for all $\tau \in [-h_M,0]$, and where the quantities appearing in (\ref{eq: aux abstract form - aux input v}) are defined by (\ref{eq: truncated model - Y}-\ref{eq: truncated model - M}). On one hand, if $X$ is a mild solution of (\ref{eq: abstract form - DE}-\ref{eq: abstract form - IC}) with $u = KY \in \mathcal{C}^1(\mathbb{R}_+;\mathbb{R}^2)$, then the developments of Section~\ref{sec: control strat and main result} show that $Y$ satisfies the ODE (\ref{eq: spectral decomposition}-\ref{eq: spectral decomposition IC}), which provides after integration (\ref{eq: aux abstract form - aux input v}) with $v = Y$. Then (\ref{eq: def aux mild solution}) holds with $\tilde{w} = w$. On the other hand, assume that $X\in\mathcal{C}^0(\mathbb{R}_+;\mathcal{H})$ satisfies (\ref{eq: def aux mild solution}). Then, we obtain from (\ref{eq: aux abstract form - aux input v}) that $v \in \mathcal{C}^1(\mathbb{R}_+;\mathbb{R}^{2N_0})$ and $\dot{v}(t) = (A_\mathrm{cl}-M) v(t) + MY(t-h(t)) + B d_b(t) + P_d(t)$. Reproducing the developments of Section~\ref{sec: control strat and main result} regarding the derivation of the ODE (\ref{eq: spectral decomposition}-\ref{eq: spectral decomposition IC}), we obtain from (\ref{eq: def aux mild solution}) that $\dot{Y}(t) = (A-M) Y(t) + M Y(t-h(t)) + B \{ K v(t) + d_b(t) \} + P_d(t)$. Consequently we have $\dot{v}(t) - \dot{Y}(t) = (A-M) (v(t) - Y(t))$ with the initial condition $v(0) - Y(0) = Y_\Phi(0) - Y_\Phi(0) = 0$. Thus $v = Y \in \mathcal{C}^1(\mathbb{R}_+;\mathbb{R}^{2N_0})$, showing that $X$ is a mild solution of (\ref{eq: abstract form - DE}-\ref{eq: abstract form - IC}) with $u = KY$.

To conclude, it is sufficient to show the existence and uniqueness of a function $X\in\mathcal{C}^0(\mathbb{R}_+;\mathcal{H})$ satisfying (\ref{eq: def aux mild solution}-\ref{eq: aux abstract form - aux input v}). As $\Phi\in\mathcal{C}^0(\mathbb{R}_+;\mathcal{H})$ and noting that for any $k \geq 0$, $0 \leq t \leq (k+1)h_m$ implies that $- h_M < -h_m \leq t - h(t) \leq k h_m$, the existence and uniqueness of such a $X\in\mathcal{C}^0(\mathbb{R}_+;\mathcal{H})$ is immediate by a classical induction argument and~\cite[Lem.~3.1.5]{Curtain2012}. Moreover, we deduce from (\ref{eq: aux abstract form - aux input v}) that the control input is such that $u = K Y \in \mathcal{C}^1(\mathbb{R}_+;\mathbb{R}^2)$.

\section{Estimate of $\Vert M_n \Vert$}\label{annex: estimate norm Mn}

Considering $M_n$ defined by (\ref{eq: definition Mn}), we show that the following estimate holds:
\begin{equation}\label{eq: estimate norm Mn}
\forall n \geq 1 , \quad \Vert M_n \Vert \leq m_n \triangleq \dfrac{\sqrt{2} \alpha \gamma}{\sqrt{(\alpha^2-1) n^4 \pi^4 + \beta}} .
\end{equation}
To do so, we recall that, for any real matrix $P$, $\Vert P \Vert = \sqrt{\lambda_\mathrm{M}(P^\top P)}$. Let $a \neq 0$ and consider the matrix
\begin{equation*}
P_a = \begin{bmatrix} 1 & a \\ -1/a & -1 \end{bmatrix} .
\end{equation*}
As the eigenvalues of $P_a^\top P_a$ are given by 0 and $2+a^2+1/a^2$, we infer $\Vert P_a \Vert = \sqrt{2+a^2+1/a^2}$. Thus we have
\begin{equation*}
\Vert M_n \Vert = \dfrac{\gamma}{\vert \lambda_{n,-1} - \lambda_{n,+1} \vert} \sqrt{2 + \dfrac{k_{n,-1}^2}{k_{n,+1}^2} + \dfrac{k_{n,+1}^2}{k_{n,-1}^2}} .
\end{equation*} 
As $\alpha > 1$, straightforward computations show that
\begin{align*}
2 + \dfrac{k_{n,-1}^2}{k_{n,+1}^2} + \dfrac{k_{n,+1}^2}{k_{n,-1}^2}
& = 4 \alpha^2 \dfrac{(\alpha n^4 \pi^4 + \beta/(2\alpha))^2}{\alpha^2 n^8 \pi^8 + \beta^2/4} \\
& \leq 4 \alpha^2 \dfrac{(\alpha n^4 \pi^4 + \beta/2)^2}{\alpha^2 n^8 \pi^8 + \beta^2/4} \leq 8 \alpha^2 ,
\end{align*}
where is has been used that $(a+b)^2 \leq 2(a^2+b^2)$ for all $a,b\in\mathbb{R}$. The claimed estimate (\ref{eq: estimate norm Mn}) follows from (\ref{eq: eigenvalue U0}).

\section{Establishment of estimate (\ref{eq: rough estimate Sc over [0,hm]})}\label{annex: rough estimate Sc over [0,hm]}

In this appendix, we always consider integers $n \geq N_0 + 1$. Thus we have $\lambda_{n,\epsilon} \leq \lambda_{N_0 +1,+1} = - 2 \eta < 0$ for all $\epsilon\in\{-1,+1\}$. From the definition of $v_n$ given by (\ref{eq: def vn(t)}) and the ODE (\ref{eq: modified spectral decomposition for stab analysis inf dim residual dynamics}), we have, for all $t \geq 0$,
\begin{align*}
c_n(t) & = e^{\Lambda_n t} c_n(0) - \int_0^t e^{\Lambda_n (t-\tau)} M_n c_n(\tau) \diff\tau \\
& \phantom{=}\; + \int_0^t e^{\Lambda_n(t-\tau)} M_n \tilde{c}_n(\tau) \diff\tau \\
& \phantom{=}\; + \int_0^t e^{\Lambda_n(t-\tau)} \{q_n(\tau) - \Lambda_n r_n(\tau) + p_{d,n}(\tau) \} \diff\tau ,
\end{align*}
where $\tilde{c}_n(t) = c_n(t-h(t))$. Noting that, for all $t \geq 0$, $\Vert e^{\Lambda_n t} \Vert = e^{\lambda_{n,+1} t} \leq e^{-2\eta t} \leq 1$, we have 
\begin{align*}
& \Vert c_n(t) \Vert \\
& \leq \Vert c_n(0) \Vert + m_{N_0+1} \mathcal{I}(\Vert c_n \Vert,0,t) \\
& \phantom{=}\;  + m_{N_0+1} \mathcal{I}(\Vert \tilde{c}_n \Vert,0,t) + \left\Vert \Lambda_n \int_0^t e^{\Lambda_n (t-\tau)} r_n(\tau) \diff\tau \right\Vert \\
& \phantom{=}\; + \mathcal{I}(\Vert q_n \Vert,0,t) + \mathcal{I}(\Vert p_{d,n} \Vert,0,t) 
\end{align*}
for all $t \geq 0$, where we have used (\ref{eq: estimate norm Mn}) and with the notation $\mathcal{I}$ defined by (\ref{eq: def integral I}). Then, we obtain for all $t \geq 0$
\begin{align*}
\Vert c_n(t) \Vert^2
& \leq 6 \Vert c_n(0) \Vert^2 + \dfrac{3 m_{N_0+1}^2}{\eta} \mathcal{I}(\Vert c_n \Vert^2,0,t) \\
& \phantom{\leq}\; + \dfrac{3 m_{N_0+1}^2}{\eta} \mathcal{I}(\Vert \tilde{c}_n \Vert^2,0,t) \\
& \phantom{\leq}\; + 6 \left\Vert \Lambda_n \int_0^t e^{\Lambda_n (t-\tau)} r_n(\tau) \diff\tau \right\Vert^2 \\
& \phantom{\leq}\; + \dfrac{3}{\eta} \mathcal{I}(\Vert q_n \Vert^2,0,t) + \dfrac{3}{\eta} \mathcal{I}(\Vert p_{d,n} \Vert^2,0,t)  ,
\end{align*}
where estimate (\ref{eq: integral I - estimate 2}) has been used. Summing for $n \geq N_0 + 1$ and then using (\ref{eq: integral I - estimate 1}) with $\kappa = 0$, we have, for all $t \geq 0$,
\begin{align*}
S_c(t)
& \leq 6 S_c(0) + \dfrac{3 m_{N_0+1}^2}{\eta} \mathcal{I}(S_c,0,t) \\
& \phantom{\leq}\; + \dfrac{3 m_{N_0+1}^2}{\eta} \mathcal{I}( S_c(\cdot - h(\cdot) ),0,t) \\
& \phantom{\leq}\; + 6 \sum\limits_{n \geq N_0 + 1} \left\Vert \Lambda_n \int_0^t e^{\Lambda_n (t-\tau)} r_n(\tau) \diff\tau \right\Vert^2 \\
& \phantom{\leq}\; + \dfrac{3}{\eta} \mathcal{I}(S_q,0,t) + \dfrac{3}{\eta} \mathcal{I}(S_p,0,t) \\
& \leq 6 S_c(0) + \xi_2 \sup\limits_{\tau\in[0,t]}S_c(\tau) + \xi_2 \sup\limits_{\tau\in[-h_M,t-h_m]}S_c(\tau) \\
& \phantom{\leq}\; + 6 \sum\limits_{n \geq N_0 + 1} \left\Vert \Lambda_n \int_0^t e^{\Lambda_n (t-\tau)} r_n(\tau) \diff\tau \right\Vert^2 \\
& \phantom{\leq}\; 
+ \xi_3 \sup\limits_{\tau\in[0,t]}S_q(\tau) + \xi_3 \sup\limits_{\tau\in[0,t]}S_p(\tau) 
\end{align*}
with $\xi_2 = 3 m_{N_0+1}^2/(2\eta^2)$ and $\xi_3 = 3/(2\eta^2)$. Using assumption (\ref{eq: thm infinite dim residual - small gain condition}), we infer that $\xi_2 < 1$ and thus
\begin{align*}
\sup\limits_{\tau\in[0,t]}S_c(\tau)
& \leq \dfrac{6}{1-\xi_2} S_c(0) + \dfrac{\xi_2}{1-\xi_2} \sup\limits_{\tau\in[-h_M,t-h_m]}S_c(\tau) \\
& \phantom{\leq}\; + \dfrac{6}{1-\xi_2} \sum\limits_{n \geq N_0 + 1} \left\Vert \Lambda_n \int_0^t e^{\Lambda_n (t-\tau)} r_n(\tau) \diff\tau \right\Vert^2 \\
& \phantom{\leq}\; + \dfrac{\xi_3}{1-\xi_2} \sup\limits_{\tau\in[0,t]}S_q(\tau) + \dfrac{\xi_3}{1-\xi_2} \sup\limits_{\tau\in[0,t]}S_p(\tau) 
\end{align*}
for all $t \geq 0$. Recalling that $S_c(0) \leq \Vert \Phi(0) \Vert^2 / m_R$ and $S_p(t) \leq \Vert d_d(t) \Vert^2 / m_R$, the use of estimates (\ref{eq: estimate Q}) and (\ref{eq: estimate series integral rn}) for $t_0 = 0$ show that 
\begin{align}
\sup\limits_{\tau\in[0,t]}S_c(\tau)
& \leq \xi_4 \sup\limits_{\tau\in[-h_M,t-h_m]} S_c(\tau)
+ \xi_5 \Vert\Phi\Vert_{1,h_M}^2 \label{eq: annex rough estimate Sc} \\
& \phantom{\leq}\; + \xi_6 \sup\limits_{\tau\in[0,t]} \Vert d_d(\tau) \Vert 
+ \xi_7 \sup\limits_{\tau\in[0,t]} \Vert d_b(\tau) \Vert , \nonumber
\end{align}
for all $t \geq 0$ with $\xi_4 = \xi_2 / (1-\xi_2)$, $\xi_5 = (6/m_R + 6 \gamma_{3,1} + \xi_3 \gamma_{1,1}) / (1-\xi_2)$, $\xi_6 = (6 \gamma_{3,2} + \xi_3 \gamma_{1,2} + \xi_3/m_R) / (1-\xi_2)$, and $\xi_7 = (6 \gamma_{3,3} + \xi_3 \gamma_{1,3}) / (1-\xi_2)$. Now, as for any $\tau \in [-h_m,0]$, $\Phi(\tau) = \Phi(0) + \int_0^\tau \dot{\Phi}(s) \diff s$, we infer
\begin{align*}
\Vert \Phi(\tau) \Vert & \leq \Vert \Phi(0) \Vert + \int_{-h_M}^0 \Vert \dot{\Phi}(s) \Vert \diff s \\
& \leq \Vert \Phi(0) \Vert + \sqrt{h_M} \sqrt{\int_{-h_M}^0 \Vert \dot{\Phi}(s) \Vert^2 \diff s} \\
& \leq \left( 1 + \sqrt{h_M} \right) \Vert \Phi \Vert_{1,h_M} ,
\end{align*}
and thus we have the estimate $S_c(\tau) \leq \Vert \Phi(\tau) \Vert^2 / m_R  \leq \left( 1 + \sqrt{h_M} \right)^2 \Vert \Phi \Vert_{1,h_M}^2 / m_R$ for all $\tau \in [-h_M,0]$. Based on (\ref{eq: annex rough estimate Sc}), one can recursively estimate the quantities $\sup\limits_{\tau\in[0,k h_m]} S_c(\tau)$ for increasing values of $k \in \mathbb{N}^*$. This procedure provides the claimed estimate (\ref{eq: rough estimate Sc over [0,hm]}) when reaching $k = \lceil h_M/h_m \rceil \geq 1$.


\ifCLASSOPTIONcaptionsoff
  \newpage
\fi



\bibliographystyle{IEEEtranS}
\nocite{*}
\bibliography{IEEEabrv,mybibfile}

\end{document}